\documentclass[12pt]{article}
\usepackage{amsmath}
\usepackage{amsthm}
\usepackage{amsfonts}
\usepackage{amssymb}
\usepackage[all]{xypic}

\makeatletter

\newdimen\p@renwd
\newdimen\tmp
\setbox0=\hbox{B} \p@renwd=\wd0 % width of the big left (
\def\bm#1{\begingroup \m@th
  \setbox\z@\vbox{\def\cr{\crcr\noalign{\kern2\p@\global\let\cr\endline}}%
    \ialign{$##$\hfil\kern2\p@\kern\p@renwd&\thinspace\hfil$##$\hfil
      &&\quad\hfil$##$\hfil\crcr
      \omit\strut\hfil\crcr\noalign{\kern-\baselineskip}%
      #1\crcr\omit\strut\cr}}%
  \setbox\tw@\vbox{\unvcopy\z@\global\setbox\@ne\lastbox}%
  \setbox\tw@\hbox{\unhbox\@ne\unskip\global\setbox\@ne\lastbox}%
  \setbox\tw@\hbox{$\kern\wd\@ne\kern-\p@renwd\left(\kern-\wd\@ne
    \global\setbox\@ne\vbox{\box\@ne\kern2\p@}%
    \vcenter{\tmp=-\ht\@ne\unvbox\z@\kern-\baselineskip\kern\tmp}\,\right)$}%
  \null\;\vbox{\kern\ht\@ne\box\tw@}\endgroup}

\makeatother

\newcounter{appendix}
\def\appendix#1{\stepcounter{appendix}\goodbreak\vspace{12pt plus 3pt}
{\begin{center} 
\Huge Appendix \Alph{appendix}\par\vspace{8pt plus 2pt}
\bfseries\Large{#1}\end{center}}
\addcontentsline{toc}{section}{Appendix \Alph{appendix}. #1}\setcounter{subsection}{0}%
\vspace{8pt plus 2pt}
\nobreak}

\let\phi\varphi
\let\epsilon\varepsilon
\def\C{\mathbb C}
\def\CP{\mathbb {CP}}
\def\Z{\mathbb Z}
\def\N{\mathbb N}
\def\g{\mathfrak g}
\def\e{e}
\def\F{\mathcal F}
\def\M{\mathcal M}
\def\E{\mathcal E}
\def\slg{\mathfrak{sl}}

\def\mmod{\ \mathrm{mod}\ }
\let\wt\widetilde
\let\wh\widehat
\def\Im{\mathop{\mathrm{Im}}}
\def\sn{\mathop{\mathrm{sn}}}
\def\dn{\mathop{\mathrm{dn}}}
\def\cn{\mathop{\mathrm{cn}}}
\def\reg{\mathop{\mathrm{reg}}}
\def\Aut{\mathop{\mathrm{Aut}}}
\def\Mat{\mathop{\mathrm{Mat}}\nolimits}
\def\Dr{\mathop{\mathrm{Dr}}\nolimits}
\def\End{\mathop{\mathrm{End}}}
\def\NOD{\mathop{\mathrm{gcd}}}

\newtheorem{prop}{Proposition}

\theoremstyle{remark}

\newtheorem*{remark}{Remark}
\newtheorem*{remarks}{Remarks}

\theoremstyle{definition}

\newtheorem*{definition}{Definition}

\title{Elliptic algebras}
\author{Alexander Odesskii}

\date{}

\begin{document}

\sloppy

\maketitle

\begin{abstract}
The survey is devoted to associative $\Z_{\ge0}$-graded algebras presented by $n$
generators and $\frac{n(n-1)}2$ quadratic relations and satisfying the so-called 
Poincare-Birkhoff-Witt condition (PBW-algebras). We consider examples of such 
algebras depending on two continuous parameters (namely, on an elliptic curve 
and a point on this curve) which are flat deformations of the polynomial ring in 
$n$ variables. Diverse properties of these algebras are described, together with
their relations to integrable systems, deformation quantization, moduli spaces 
and other directions of modern investigations.
\end{abstract}

\newpage

\tableofcontents

\section*{Introduction}
\addcontentsline{toc}{section}{Introduction}

In the paper [45] devoted to study the $XYZ$-model and the representations of 
the corresponding algebra of monodromy matrices, Sklyanin introduced the family 
of associative algebras with four generators and six quadratic relations which 
are nowadays called Sklyanin algebras (see also Appendix D.1). The algebras of 
this family are naturally indexed by two continuous parameters, namely, by an 
elliptic curve and a point on this curve, and each of them is a flat deformation 
of the polynomial ring in four variables in the class of $\Z_{\ge0}$-graded 
associative algebras. On the other hand, a family of algebras with three generators 
(and three quadratic relations) with the same properties arose in [2], [34] (see also
[52]). In what follows it turned out (see [10], [17]-[22], [32]-[38]) that such 
algebras exist for arbitrarily many generators. The algebras in question are 
associative algebras of the following form. Let $V$ be a linear space of dimension 
$n$ over the field $\C$. Let $L\subset  V\otimes  V$ be a subspace of dimension 
$\frac{n(n-1)}2$. Let us construct an algebra $A$ with the space of generators 
$V$ and the space of defining relations $L$, that is, $A=T^*V/(L)$, where $T^*V$ 
is the tensor algebra of the space $V$ and $(L)$ is the two-sided ideal generated 
by $L$. It is clear that the algebra $A$ is $\Z_{\ge0}$-graded because the ideal 
$(L)$ is homogeneous. We have $A=\C\oplus  A_1\oplus A_2\oplus\dots$, where 
$A_1=V$,  $A_2=V\otimes
V/L$, $A_3=V\otimes V\otimes V/V\otimes L+L\otimes V$, etc.

\begin{definition}
We say that $A$ is a PBW-algebra (or satisfies the Poincare-Birkhoff-Witt condition) 
if $\dim A_\alpha=\frac{n(n+1)\dots(n+\alpha-1)}{\alpha!}$.
\end{definition}

Thus, a PBW-algebra is an algebra with $n$ generators and $\frac{n(n-1)}2$ 
quadratic relations for which the dimensions of the graded components are equal to 
those of the polynomial ring in $n$ variables.

Algebras of this kind arise in diverse areas of mathematics: in the theory of 
integrable systems [45], [46], [28], [9], moduli spaces [20], deformation quantization 
[12], [26], non-commutative geometry [2], [3], [11], [27], [47]-[49], [51], 
cohomology of algebras [8], [29], [41]-[44], [50], and quantum groups and 
$R$-matrices [45], [46], [25], [16], [14], [23], [31]. See Appendix D.

Since there are no classification results in the theory of PBW-algebras (for $n>3$), 
we deal with specific examples only. The known examples can conditionally be 
divided into two classes, namely, rational and elliptic algebras. Let us present 
examples of rational algebras.

1. {\it  Skew polynomials.} This is the algebra with the generators $\{x_i;i=1,\dots,n\}$ 
and the relations $x_ix_j=q_{i,j}x_jx_i$, where $i<j$ and $q_{i,j}\ne0$.

One can readily see that the monomials 
$\{x_1^{\alpha_1}\dots x_n^{\alpha_n};\alpha_1,\dots,\alpha_n\in\Z_{\ge0}\}$ 
form a basis of the algebra of skew polynomials, which implies the PBW condition. 
Since $q_{i,j}$ are arbitrary non-zero numbers, we have obtained an 
$\frac{n(n-1)}2$-parameter family of algebras.

2. {\it Projectivization of Lie algebras.} Let $\g$ be a Lie algebra of dimension $n-1$ 
with a basis $\{x_1,\dots,x_{n-1}\}$. We construct an algebra with $n$ generators 
$\{c,x_1,\dots,x_{n-1}\}$ and the relations $cx_i=x_ic$ and 
$x_ix_j-x_jx_i=c[x_i,x_j]$.

The condition PBW follows from the Poincare-Birkhoff-Witt theorem for the algebra 
$\g$. 

3. {\it Drinfeld algebra.} A new realization of the quantum current algebra 
$U_q(\wh{\slg}_2)$ was suggested in [13] (see also [25]). Namely, the generators 
$x_k^\pm,h_k$ ($k\in\Z$) similar to the ordinary basis of the Lie algebra 
$\wh{\slg}_2$ were introduced. It is assumed that the elements $x_k^+$ satisfy 
the quadratic relations 
\begin{equation}
x_{k+1}^+x_l^+-q^2x_l^+x_{k+1}^+=q^2x_k^+x_{l+1}^+-x_{l+1}^+x_k^+.
\end{equation}
The elements $x_k^-$ satisfy similar relations. The algebra 
$\Dr_n(q)\subset U_q(\wh{\slg}_2)$ generated by $x_1^+,\dots,x_n^+$,
$n\in\N$, $q\in\C^*$, is a PBW-algebra.

In the elliptic case the algebra depends on two continuous parameters, namely, an 
elliptic curve $\E$ and a point $\eta\in\E$. Just these algebras are the subject 
of our survey. Their structure constants are elliptic functions of $\eta$ with 
modular parameter $\tau$. Our main example is given by the algebras 
$Q_{n,k}(\E,\eta)$, where $n\ge3$ is the number of generators, $k$ is a positive 
integer coprime to $n$, and $1\le  k<n$. We define the algebra $Q_{n,k}(\E,\eta)$ 
by the generators $\{x_i;i\in\Z/n\Z\}$ and the relations 
\begin{equation}
\sum_{r\in\Z/n\Z}\frac{\theta_{j-i+r(k-1)}(0)}
{\theta_{kr}(\eta)\theta_{j-i-r}(-\eta)}x_{j-r}x_{i+r}=0.
\end{equation}
The structure of these algebras depends on the expansion of the number $n/k$ in the 
continued fraction, and therefore we first study the simplest case $k=1$ and then 
pass to the general case. The fact that the algebra $Q_{n,k}(\E,\eta)$ belongs to 
the class of PBW-algebras is proved only for generic parameters $\E$ and $\eta$ 
(see \S2.6 and \S3). However, we conjecture that this holds for any $\E$ and $\eta$. 
A possible way to prove this conjecture is to produce an analog of the functional 
realization (see \S2.1) for arbitrary $k$ by using the constructions in \S5.

As we consider, the algebras $Q_{n,k}(\E,\eta)$ are a typical example of elliptic 
algebras; however, they are far from exhausting the list of all elliptic algebras. 
The simplest example of an elliptic algebra that does not belong to this class 
(and even is not a deformation of the polynomial ring) can be constructed as follows. 
Let the group $(\Z/2\Z)^2$ with the generators $g_1,g_2$ act by automorphisms on 
the algebra $Q_4(\E,\eta)$ as follows:
$g_1(x_i)=x_{i+2}$,
$g_2(x_i)=(-1)^ix_i$.
The same group acts on the algebra of ($2\times2$) matrices,
$g_1(\gamma)=\begin{pmatrix}-1&0\\0&1\end{pmatrix}\gamma\begin{pmatrix}-1&0\\0&1\end{pmatrix}^{-1}$, 
$g_2(\gamma)=\begin{pmatrix}0&1\\1&0\end{pmatrix}\gamma\begin{pmatrix}0&1\\1&0\end{pmatrix}^{-1}$.
This gives an action on the tensor product of associative algebras
$Q_4(\E,\eta)\otimes\Mat_2$. Let $\wt Q_4(\E,\eta)\subset
Q_4(\E,\eta)\otimes\Mat_2$ consist of elements invariant with respect to the group
action. One can readily see that the dimension of the graded components of 
$\wt Q_4(\E,\eta)$ coincide with those of $Q_4(\E,\eta)$, and therefore 
$\wt Q_4(\E,\eta)$ is a PBW-algebra. For another example of PBW-algebra (with $3$ 
generators), see the end of \S1.

Let us now describe one of the main constructions of PBW-algebras. Let 
$\lambda(x,y)$ be a meromorphic function of two variables. We construct an 
associative graded algebra $\F_\lambda$ as follows. Let the underlying linear space 
of $\F_\lambda$ coincide with $\F_\lambda=\C\oplus   F_1\oplus  F_2\oplus\dots$, 
where $F_1=\{f(u)\}$ is the space of meromorphic functions of one variable and 
$F_\alpha=\{f(u_1,\dots,u_\alpha)\}$ is the space of symmetric meromorphic functions 
of $\alpha$ variables. The space $F_\alpha$ is a natural extension of the symmetric 
power  $S^\alpha F_1$. The multiplication in the algebra $\F_\lambda$ is defined as 
follows: for $f\in F_\alpha$, and $g\in F_\beta$ the product $f*g\in
F_{\alpha+\beta}$ is 

\medskip
$f*g(u_1,\dots,u_{\alpha+\beta})=$
\begin{equation}
=\frac1{\alpha!\beta!}\sum_{\sigma\in
S_{\alpha+\beta}}f(u_{\sigma_1},\dots,u_{\sigma_\alpha})
g(u_{\sigma_{\alpha+1}},\dots,u_{\sigma_{\alpha+\beta}})
\prod_{\begin{subarray}{c}1\le                i\le\alpha\\\alpha+1\le
j\le\alpha+\beta\end{subarray}}\lambda(u_{\sigma_i},u_{\sigma_j}).
\end{equation}
In particular, if $f,g\in F_1$, then 
\begin{equation}
f*g(u_1,u_2)=f(u_1)g(u_2)\lambda(u_1,u_2)+
f(u_2)g(u_1)\lambda(u_2,u_1).
\end{equation}
One can readily see that the multiplication $*$ is associative for any $\lambda(x,y)$.

We now assume that $\lambda(x,y)=\frac{x-qy}{x-y}$, where $q\in\C^*$. Let 
$F_1^{(n)}=\{1,u,\dots,u^{n-1}\}\subset    F_1$  be the space of polynomials 
of degree less than $n$. Let $F_\alpha^{(n)}=S^\alpha
F_1^{(n)}\subset  F_\alpha$ be the space of symmetric polynomials in $\alpha$ 
variables of degree less than $n$ with respect to any variable. One can readily see 
that $F_\alpha^{(n)}*F_\beta^{(n)}\subseteq
F_{\alpha+\beta}^{(n)}$. Therefore, the algebra $\F_\lambda^{(n)}=\oplus_\alpha
F_\alpha^{(n)}$ is a subalgebra of $\F_\lambda$. Moreover, for $q=1$ the algebra 
$\F_\lambda^{(n)}$ is the polynomial ring $S^*F_1^{(n)}$ because $\lambda(x,y)=1$ 
in this case. Therefore, the algebra $\F_\lambda^{(n)}$ is a PBW-algebra for generic 
$q$. This algebra is isomorphic to the Drinfeld algebra $\Dr_n(q)$, and an 
isomorphism is given by the rule $u^k\mapsto x_{k+1}^+$. The algebra 
$Q_n(\E,\eta)$ can be obtained in a similar way with the only modification that the 
polynomials are replaced by theta functions (see \S2.1). A similar construction [38], 
[22] enables one to construct quantum moduli spaces $\M(\E,B)$ (see Appendix D.3) 
for any Borel subgroup $B$. The construction of algebras $Q_{n,k}(\E,\eta)$ for 
$k>1$ (and, more generally, quantum moduli spaces $\M(\E,P)$ for a parabolic 
subgroup $P$) is more complicated and involves exchange algebras (see \S5 and [21]) 
or elliptic $R$-matrices (see \S4).

Let us now describe the contents of the survey. In \S1 we describe the simplest 
elliptic PBW-algebras, namely, algebras $Q_3(\E,\eta)$ with three generators. These 
algebras were studied in many papers, see, for instance, [2], [3]. The section is 
of illustrative nature; we intend to explain some methods of studying elliptic 
algebras by the simplest example. The main attention in the survey is paid to the 
algebras $Q_n(\E,\eta)$, which are discussed in \S2. We give an explicit construction 
of these algebras, present natural families of their representations (which are 
studied in [19] in more detail), and describe the symplectic leaves of the 
corresponding Poisson algebra (we recall that $Q_n(\E,0)$ is the polynomial ring 
in $n$ variables). 

The structure of the algebras $Q_{n,k}(\E,\eta)$, $k>1$, is more complicated, and 
the detailed description of their properties is beyond the framework of the survey 
(see [35], [20]). The main properties of these algebras are described in \S3. In 
\S4 we explain the relationship between these algebras and Belavin's elliptic 
$R$-matrices. In \S5 we establish a relation of the algebras $Q_{n,k}(\E,\eta)$ to 
the so-called exchange algebras (see (24), (25), and also [36], [24], [33]). In 
Appendices A, B, C we present the notation we need and the properties of theta 
functions of one and several variables. Appendix D contains a brief survey of 
relations of elliptic algebras with other areas of mathematics. We tried to make 
this part independent of the main text.

In conclusion we say a few words concerning the facts that remain outside the survey 
but are immediately connected with its topic. In [37] the algebras 
$Q_{n,k}(\E,\eta)$ are studied provided that $\eta\in\E$ is a point of finite order. 
In this case the properties of the algebras $Q_{n,k}(\E,\eta)$ are similar to those 
of quantum groups if $q$ is a root of unity; in particular, these algebras are 
finite-dimensional over the centre. In [32] we study rational degenerations of the 
algebras $Q_{n,k}(\E,\eta)$ occurring if the elliptic curve $\E$ degenerates into 
the union of several copies of $\CP^1$ or into $\CP^1$ with a double point.

The algebras $Q_{n,k}(\E,\eta)$ are obtained when quantizing the components of the 
moduli spaces $\M(P,\E)$ (see Appendix D.3) that are isomorphic to $\CP^{n-1}$. The 
quantization of other components leads to elliptic algebras of more general form. 
These algebras were constructed in [38], [22] if $P$ is a Borel subgroup of an 
arbitrary group $G$. The case in which $P\subset   GL_m$ is an arbitrary parabolic 
subgroup of $GL_m$ is studied in [21].

The symplectic leaves of a Poisson manifold corresponding to the family of algebras 
$Q_{n,k}(\E,\eta)$ in a neighbourhood of $\eta=0$ and for a fixed elliptic curve 
$\E$ were studied in [20].

The corresponding Poisson algebras belong to the class of algebras with regular 
structure of symplectic leaves; these algebras were studied in [39].

\section{Algebras with three generators}

In this section we consider the simplest examples of elliptic PBW-algebras, namely, 
the algebras with three generators. Let us first study the quadratic Poisson 
structures on $\C^3$. Let $x_0,x_1,x_2$ be the coordinates on $\C^3$ and let there 
be a Poisson structure that is quadratic in these coordinates. We construct the 
polynomial $C=x_0\{x_1,x_2\}+x_1\{x_2,x_0\}+x_2\{x_0,x_1\}$. This is a homogeneous 
polynomial of degree three because the Poisson structure is quadratic. It is clear 
that the form of this polynomial is preserved under linear changes of coordinates 
(up to proportionality). Let us restrict ourselves to the non-degenerate case in 
which the equation  $C=0$ defines a non-singular projective manifold. It is clear 
that this is an elliptic curve. Moreover, by a linear change of variables one can 
reduce the polynomial $C$ to the form $C=x_0^3+x_1^3+x_2^3+3kx_0x_1x_2$, where 
$k\in\C$. In this case, as one can readily see by using the definition of $C$ and 
the Jacoby identity, the Poisson structure must be of the form (up to proportionality):
\begin{equation}
\{x_0x_1\}=x_2^2+kx_0x_1,\quad\{x_1x_2\}=x_0^2+kx_1x_2,\quad
\{x_2x_0\}=x_1^2+kx_2x_0.
\end{equation}
Moreover, $\{x_i,C\}=0$, and every central element is a polynomial in $C$. We recall 
that each Poisson manifold can be partitioned into the so-called symplectic leaves, 
which are Poisson submanifolds, and the restrictions of the Poisson structure to 
these submanifolds are non-degenerate. In our case, the symplectic leaves are as 
follows:

1) the origin $x_0=x_1=x_2=0$;

2) the homogeneous manifold $C=0$ without the origin;

3) the manifolds $C=\lambda$, where $\lambda\in\C$, $\lambda\ne0$.

It is clear that our Poisson structure admits the automorphisms 
$x_i\mapsto\epsilon^ix_i$ and $x_i\mapsto  x_{i+1}$, where $\epsilon^3=1$,
$i\in\Z/3\Z$. It is natural to assume that the quantization of the Poisson structure 
(see Appendix D.2) is the family of associative algebras with the generators 
$x_0,x_1,x_2$ and three quadratic relations admitting the same automorphisms. 
However, each generic three-dimensional space of quadratic relations which is 
invariant with respect to these automorphisms is of the form
\begin{equation}
\begin{aligned}
x_0x_1-qx_1x_0&=px_2^2,\\
x_1x_2-qx_2x_1&=px_0^2,\\
x_2x_0-qx_0x_2&=px_1^2,
\end{aligned}
\end{equation}
where $p,q\in\C$ are complex numbers. We denote by $A_{p,q}$ the algebra with the 
generators $x_0,x_1,x_2$ and the defining relations (6). It is clear that the 
algebra $A_{p,q}$ is $\Z_{\ge0}$-graded, that is, $A_{p,q}=\C\oplus    F_1\oplus
F_2\oplus\dots$, where $F_\alpha F_\beta\subseteq F_{\alpha+\beta}$. Here $F_\alpha$ 
stands for the linear space spanned by the (non-commutative) monomials in 
$x_0,x_1,x_2$ of degree $\alpha$. It is natural to expect that the dimension of 
$F_\alpha$ is equal to that of the space of polynomials in three variables of degree 
$\alpha$, that is, $\dim F_\alpha=\frac{(\alpha+1)(\alpha+2)}2$.

Moreover, the Poisson algebra (5) has a central function 
$C=x_0^3+x_1^3+x_2^3+3kx_0x_1x_2$, and the centre is generated by the element $C$. 
Therefore, it is natural to expect that for generic $p$ and $q$ the algebra 
$A_{p,q}$ has a central element of the form $C_{p,q}=\phi x_0^3+\psi x_1^3+\mu
x_2^3+\lambda x_0x_1x_2$, where $\phi,\psi,\mu$, and $\lambda$ are functions of $p$ 
and $q$ (one can verify the existence of an element $C_{p,q}$ by the immediate 
calculation), and the centre is generated by $C_{p,q}$.

The standard technique of proving such assertions (for instance, the Poincare-Birkhoff
-Witt theorem for Lie algebras) makes use of the filtration on an algebra and the 
study of the graded adjoint algebra. In our case the algebra is already graded, and 
one cannot proceed by the ordinary induction on the terms of lesser filtration; 
therefore we use another technique. Namely, we shall study a certain class of modules 
over the algebra $A_{p,q}$ and try to obtain results on the algebra $A_{p,q}$ by 
using an information on the modules. The following class of modules is useful for 
our purposes.

\begin{definition}
A module over a $\Z_{\ge0}$-graded algebra $A$ is said to be linear if it is 
$\Z_{\ge0}$-graded as an $A$-module, generated by the space of degree $0$, and the 
dimensions of all components are equal to $1$.
\end{definition}

Let us study the linear modules over the algebra $A_{p,q}$. By definition, a linear 
module $M$ admits a basis $\{v_\alpha,\alpha\ge0\}$ with the following action of 
the generators:
$$
x_0v_\alpha=x_\alpha     v_{\alpha+1},\quad      x_1v_\alpha=y_\alpha
v_{\alpha+1},\quad x_2v_\alpha=z_\alpha v_{\alpha+1},
$$
where $x_\alpha,y_\alpha,z_\alpha$ are sequences, and $x_\alpha,y_\alpha,z_\alpha$ 
do not vanish simultaneously for any $\alpha$ (we want $M$ be generated by $v_0$). 
A change of the basis of the form $v_\alpha\to\lambda_\alpha       v_\alpha$ 
multiplies the triple $(x_\alpha,y_\alpha,z_\alpha)\in\C^3$ by 
$\frac{\lambda_{\alpha+1}}{\lambda_\alpha}$, that is, the module $M$ is defined by 
the sequence of points $(x_\alpha:y_\alpha:z_\alpha)\in\CP^2$ uniquely up to 
isomorphism of graded modules. It is clear that a sequence of points 
$(x_\alpha:y_\alpha:z_\alpha)\in\CP^2$ defines a module over the algebra $A_{p,q}$ 
if and only if the relations (6) hold for the operators on $M$ corresponding to this 
sequence. This is equivalent to the following relations:
\begin{equation}
\begin{aligned}
x_{\alpha+1}y_\alpha-qy_{\alpha+1}x_\alpha&=pz_{\alpha+1}z_\alpha,\\
y_{\alpha+1}z_\alpha-qz_{\alpha+1}y_\alpha&=px_{\alpha+1}x_\alpha,\\
z_{\alpha+1}x_\alpha-qx_{\alpha+1}z_\alpha&=py_{\alpha+1}y_\alpha.
\end{aligned}
\end{equation}

The relations (7) form a system of linear equations for $x_\alpha,y_\alpha,z_\alpha$ 
which has a non-zero solution (by the assumption on the module $M$), and therefore 
the determinant $\begin{vmatrix}-qy_{\alpha+1}&x_{\alpha+1}&-pz_{\alpha+1}\\
-px_{\alpha+1}&-qz_{\alpha+1}&y_{\alpha+1}\\
z_{\alpha+1}&-py_{\alpha+1}&-qx_{\alpha+1}\end{vmatrix}$ must vanish. Similarly, the 
relations (7) form a system of linear equations on 
$x_{\alpha+1},y_{\alpha+1},z_{\alpha+1}$ that has a non-zero solution, and therefore 
$\begin{vmatrix}y_\alpha&-qx_\alpha&-pz_\alpha\\
-px_\alpha&z_\alpha&-qy_\alpha\\
-qz_\alpha&-py_\alpha&x_\alpha
\end{vmatrix}=0$. One can readily see that these determinants give the same cubic 
polynomial in three variables. We see that for any $\alpha\ge0$ the point with 
the coordinates $(x_\alpha:y_\alpha:z_\alpha)$ belongs to the cubic 
\begin{equation}
x_\alpha^3+y_\alpha^3+z_\alpha^3+\frac{p^3+q^3-1}{pq}
x_\alpha y_\alpha z_\alpha=0.
\end{equation}
Moreover, if a point $(x_\alpha:y_\alpha:z_\alpha)$ belongs to this cubic, then, 
solving the system of linear equations (7) with respect to 
$x_{\alpha+1},y_{\alpha+1},z_{\alpha+1}$, we obtain a new point  
$(x_{\alpha+1}:y_{\alpha+1}:z_{\alpha+1})$ on the same cubic (because the 
determinant of the system (7) must be equal to $0$). Thus, the system (7) defines 
an automorphism of the projective manifold (8). Let us choose some 
$k=\frac{p^3+q^3-1}{pq}$. Then, varying $q$, we obtain a one-parameter family 
of automorphisms of the projective curve in $\CP^2$ given by the equation 
$x^3+y^3+z^3+kxyz=0$. As is known, for generic $k$ this equation defines an 
elliptic curve. Let this curve be $\E=\C/\Gamma$, where $\Gamma$ is an integral 
lattice generated by $1$ and $\tau$, where $\Im\tau>0$. The parameter $k$ is a 
function of $\tau$. If $k$ is chosen, then, passing to the limit as $q\to1$, we 
see that $p\to0$, and the automorphism defined by (7) tends to the identity 
automorphism. Therefore, our family of automorphisms of the elliptic curve $\E$ 
given by the equation (8) is a deformation of the identity automorphism. Thus, 
every automorphism of this family is a translation, of the form $u\to u+\eta$, 
where $u,\eta\in\E=\C/\Gamma$. Let $u_\alpha\in\E=\C/\Gamma$ be a point with the 
coordinates $(x_\alpha:y_\alpha:z_\alpha)$. We see that $u_{\alpha+1}=u_\alpha+\eta$, 
where $\eta$ depends only on the algebra, that is, on $p$ and $q$. Hence, 
$u_\alpha=u+\alpha\eta$, where $u\in\E$ is the parameter of the module $M$. We 
have obtained the following result.

\begin{prop}
The linear modules over the algebra $A_{p,q}$ are indexed by a point of the 
elliptic curve $\E\subset\CP^2$ given by the equation $x^3+y^3+z^3+k_{p,q}xyz=0$, 
where $k_{p,q}=\frac{p^3+q^3-1}{pq}$. The module $M_u$ corresponding to a point 
$u\in\E$ is given by the formulas
$$
x_0v_\alpha=x_\alpha     v_{\alpha+1},\quad      x_1v_\alpha=y_\alpha
v_{\alpha+1},\quad x_2v_\alpha=z_\alpha v_{\alpha+1},
$$
where $(x_\alpha:y_\alpha:z_\alpha)$ are the coordinates of the point 
$u+\alpha\eta\in\E$. Here the shift $\eta$ is determined by $p$ and $q$.
\end{prop}

We note that, when studying linear modules, for an algebra $A_{p,q}$ we have 
constructed both an elliptic curve $\E\subset\CP^2$ and a point $\eta\in\E$. In 
what follows we shall see that, conversely, the algebra $A_{p,q}$ can be  
reconstructed from $\E$ and $\eta$. Thus, two continuous parameters, $\E$ 
(that is $\tau$) and $\eta$, give a natural parametrization of the algebras 
$A_{p,q}$. Therefore, we change the notation and denote the algebra $A_{p,q}$ by 
$Q_3(\E,\eta)$.

Let us now apply a uniformization of the elliptic curve $\E\subset\CP^2$ given 
by the equation (8) by theta functions of order three (see Appendix A). A point 
$u\in\E=\C/\Gamma$ has the coordinates $(\theta_0(u):\theta_1(u):\theta_2(u))\in\CP^2$. 
In this notation, the module $M_u$ is given by the formulas
$$
x_0v_\alpha=\theta_0(u+\alpha\eta)v_{\alpha+1},\quad
x_1v_\alpha=\theta_1(u+\alpha\eta)v_{\alpha+1},\quad
x_2v_\alpha=\theta_2(u+\alpha\eta)v_{\alpha+1}.
$$
Let $\e$ be the linear operator in the space with basis $\{v_\alpha,\alpha\ge0\}$ 
given by the formula $\e v_\alpha=v_{\alpha+1}$. Let $u$ be the diagonal operator 
in the same space such that $\e u=(u-\eta)\e$. We have 
$uv_\alpha=(u_0+\alpha\eta)v_\alpha$ for some $u_0\in\C$. It is clear that the  
generators of the algebra $Q_3(\E,\eta)$ in the representation $M_u$ become 
$$
x_0=\theta_0(u)\e,\quad
x_1=\theta_1(u)\e,\quad
x_2=\theta_2(u)\e.
$$

This gives the following reformulation of the description of linear modules.

\begin{prop}
Let us consider the $\Z_{\ge0}$-graded algebra $B(\eta)=\C\oplus
B_1\oplus B_2\oplus\dots$, where  $B_\alpha=\{f(u)\e^\alpha\}$, $f$ ranges over 
all holomorphic functions, and the multiplication is given by the formula \emph:
$f(u)\e^\alpha\cdot
g(u)\e^\beta=f(u)g(u-\alpha\eta)\e^{\alpha+\beta}$. Then there is an algebra 
homomorphism $\phi\colon Q_3(\E,\eta)\to B(\eta)$ such that 
$x_0\to\theta_0(u)\e$,           $x_1\to\theta_1(u)\e$,
$x_2\to\theta_2(u)\e$.
\end{prop}

Proposition 2 provides a lower bound for the dimension $\dim F_\alpha$ of the 
graded components of the algebra $Q_3(\E,\eta)$. Really, the homomorphism $\phi$ 
preserves the grading, that is, $\phi(F_\alpha)\subset B_\alpha$. We have
$$
\phi(x_{i_1}\dots
x_{i_\alpha})=\theta_{i_1}(u)\e\dots\theta_{i_\alpha}(u)\e=
\theta_{i_1}(u)\theta_{i_2}(u-\eta)\dots
\theta_{i_2}(u-(\alpha-1)\eta)\e^\alpha.
$$
Thus, $\phi(F_\alpha)$ is the linear space (of holomorphic functions) spanned by 
the functions $\{\theta_{i_1}(u),\dots,\theta_{i_2}(u-(\alpha-1)\eta)\};
i_1,\dots,i_\alpha=0,1,2$. It is clear that all these functions are theta 
functions of order $3\alpha$ and belong to the space 
$\Theta_{3\alpha,\frac{\alpha(\alpha-1)}23\eta}(\Gamma)$. One can readily prove 
that the image $\phi(F_\alpha)$ coincides with the entire space 
$\Theta_{3\alpha,\frac{\alpha(\alpha-1)}23\eta}(\Gamma)$, and hence 
$\dim\phi(F_\alpha)=3\alpha$. We have obtained the bound $\dim  F_\alpha\ge3\alpha$. 
On the other hand, we know that $\dim F_\alpha\le\frac{(\alpha+1)(\alpha+2)}2$ 
because the relations in $Q_3(\E,\eta)$ are deformations of the relations in the 
polynomial ring in three variables. We expect that the equality $\dim
F_\alpha=\frac{(\alpha+1)(\alpha+2)}2$ holds for generic $\tau$ and $\eta$. Let 
us compare these numbers:
$$
\begin{array}{c|c|c|c|c}
\alpha&1&2&3&4\\
\hline
\dim F_\alpha(\text{conjecture})&3&6&10&15\\
\hline
\dim\phi(F_\alpha)&3&6&9&12
\end{array}
$$
We see that the first discrepancy holds for $\alpha=3$; possibly $\phi$ has a 
one-dimensional kernel on the space $F_3$. It can be shown that, really, there is 
a cubic element $C\in Q_3(\E,\eta)$ such that $C\ne0$ and $\phi(C)=0$. The 
element $C$ turns out to be central, that is, $x_\alpha    C=Cx_\alpha$   for
$\alpha=0,1,2$. Passing to the limit as $\eta\to0$ (for a fixed $\tau$), we see 
that $C\to    x_0^3+x_1^3+x_2^3+kx_0x_1x_2$ because the $\theta_i(u)$s uniformize 
the elliptic curve, that is, $\theta_0^3+\theta_1^3+\theta_2^3+k\theta_0\theta_1\theta_2=0$. 
Further, if $C$ is central and is not a zero divisor (the latter obviously holds 
for generic $\tau$ and $\eta$), then every element $\ker\phi$ must be divisible 
by $C$ according to the dimensional considerations. Really, the graded linear 
space $\oplus_{\alpha\ge0}F_\alpha$ turns out to be not smaller than 
$\left(\oplus_{\alpha\ge0}\Theta_{3\alpha,\frac{\alpha(\alpha-1)}23\eta}\right)\otimes\C[C]$, where $\deg C=3$. 
One can readily see that the component of degree $\alpha$ of this tensor product 
of graded linear spaces is of dimension $\frac{(\alpha+1)(\alpha+2)}2$. However, 
we know that $\dim F_\alpha\le\frac{(\alpha+1)(\alpha+2)}2$, which implies 
$\dim  F_\alpha=\frac{(\alpha+1)(\alpha+2)}2$. We have obtained the following 
result.

\begin{prop}
For generic $\tau$  and  $\eta$ the algebra $Q_3(\E,\eta)$ has a cubic central 
element $C$. The quotient algebra $Q_3(\E,\eta)/(C)$ is isomorphic to 
$\oplus_{\alpha\ge0}\Theta_{3\alpha,\frac{\alpha(\alpha-1)}23\eta}(\Gamma)$, where the product 
of elements $f\in\Theta_{3\alpha,\frac{\alpha(\alpha-1)}23\eta}(\Gamma)$ and 
$g\in\Theta_{3\beta,\frac{\beta(\beta-1)}23\eta}(\Gamma)$ is given by the formula 
$f*g(u)=f(u)g(u-3\alpha\eta)$.
\end{prop}

It follows from our description of $Q_3(\E,\eta)/(C)$ that this algebra is 
centre-free for generic $\eta$. Therefore, the centre of the algebra $Q_3(\E,\eta)$ 
is generated by the element $C$.

Let us now find the relations in the algebra $Q_3(\E,\eta)$, that is, let us 
express $p$ and $q$ in term of $\tau$ and $\eta$. We have 
$x_ix_{i+1}-qx_{i+1}x_i-px_{i+2}^2=0$ (these are the relations in (6)). Applying 
the homomorphism $\phi$, we obtain
$$
\theta_i(u)\theta_{i+1}(u-\eta)-q\theta_{i+1}(u)\theta_i(u-\eta)-
p\theta_{i+2}(u)\theta_{i+2}(u-\eta)=0.
$$
Hence (see (28) in Appendix A), $q=-\frac{\theta_1(\eta)}{\theta_2(\eta)}$,
$p=-\frac{\theta_0(\eta)}{\theta_2(\eta)}$.

The similar investigation of the Sklyanin algebra with four generators (see 
Appendix D.1) gives the following result.

\begin{prop}
For a generic Sklyanin algebra $S$ with four generators and the relations (39) 
one can find an elliptic curve $\E=\C/\Gamma$ defined by two quadrics in $\CP^3$ 
and a point $\eta\in\E$ such that \emph:

\emph{1)} there is a graded algebra homomorphism $\phi\colon
S\to B(\eta)$;

\emph{2)} the image of this homomorphism in $B_\alpha$ is 
$\Theta_{4\alpha,\frac{\alpha(\alpha-1)}24\eta+\frac\alpha2}(\Gamma)$;

\emph{3)} the kernel of this homomorphism is generated by two quadratic elements 
$C_1$ and $C_2$.

Thus, $S/(C_1,C_2)=\oplus_{\alpha\ge0}\Theta_{4\alpha,\frac{\alpha(\alpha-1)}24\eta+\frac\alpha2}(\Gamma)$.
\end{prop}

The Sklyanin algebra $S$ can be reconstructed from $\E$ and $\eta$. Let us 
denote this algebra by $Q_4(\E,\eta)$.

The following natural question arises: Does there exist a similar algebra 
$Q_n(\E,\eta)$ for any $n$?

To answer this question, the information concerning linear modules is insufficient 
because these modules are too small to reconstruct the algebra $Q_n(\E,\eta)$ 
for any $n$. Really, the algebra  $Q_n(\E,\eta)$ must have the functional 
dimension $n$, whereas the linear modules are of dimension one. Therefore, these 
modules can be used only when reconstructing a quotient algebra of $Q_n(\E,\eta)$. 
To overcome these difficulties, it is natural to study more general modules. 
Namely, let us study modules over the algebra $Q_3(\E,\eta)$ with a basis 
$\{v_{i,j};i,j\in\Z_{\ge0}\}$ and such that the generators of the algebra 
$Q_3(\E,\eta)$ take any element $v_{ij}$ to a linear combination of 
$v_{i+1,j}$ and $v_{i,j+1}$. Calculations show that every such module is of 
the form
$$
x_iv_{\alpha,\beta}=\frac{\theta_i(u_1+(\alpha-2\beta)\eta)}
{\theta(u_1-u_2+3(\alpha-\beta)\eta)}v_{\alpha+1,\beta}+
\frac{\theta_i(u_2+(\beta-2\alpha)\eta)}
{\theta(u_2-u_1+3(\beta-\alpha)\eta)}v_{\alpha,\beta+1},
$$
where $i\in\Z/3\Z$, $\alpha,\beta\in\Z_{\ge0}$, and $u_1,u_2\in\C$. Thus, the 
modules of this kind are indexed by a pair of points $u_1,u_2\in\E$. If we now 
assume that the algebra $Q_n(\E,\eta)$ has analogous modules (see (15)), then 
the above information uniquely defines the algebra $Q_n(\E,\eta)$.

\begin{remarks}
1. One can pose the following more general problem. Let $M\subset\CP^{n-1}$ be 
a projective manifold and let $T$ be an automorphism of $M$. For a point $u\in  M$ 
we denote by $z_i(u)$ (where $i=0,\dots,n-1$) the homogeneous coordinates of $u$. 
Does there exist a PBW-algebra with $n$ generators $\{x_i,i=0,\dots,n-1\}$ that 
has a linear module $L_u$ (for any point $u\in   M$) given by the formula 
$x_iv_\alpha=z_i(T^\alpha
u)v_{\alpha+1}$? Here $T^\alpha u$ stands for $T(T(\dots  T(u)\dots)$. The algebras 
$Q_{n,k}(\E,\eta)$ are a solution of this problem for some $M$ and $T$, namely, 
if $M=\E^p$ is a power of a curve $\E$ and $T$ a translation (see \S5, Proposition 12). 
Here $p$ stands for the length of the expansion of $n/k=n_1-\frac1{n_2-\ldots-\frac1{n_p}}$ 
in the continued fraction.

2. Let \footnote{This example was communicated to the author by Oleg Ogievetsky 
[1], [15], [40].} $A_3$ be the algebra with the generators $x,y,z$ and the 
relations $\epsilon 
zx+\epsilon^5y^2+xz=0$, $\epsilon^2z^2+yx+\epsilon^4xy=0$, and 
$zy+\epsilon^7yz+\epsilon^8x^2=0$, where $\epsilon^9=1$. This PBW-algebra 
corresponds to the case in which $M\subset\CP^2$ is an elliptic curve given by 
the equation $x^3+y^3+z^3=0$ and $T$ is an automorphism corresponding to the 
complex multiplication on $M$. The algebra $A_3$ is not a quantization of any 
Poisson structure on $\C^3$.
\end{remarks}

\section{Algebra $Q_n(\E,\eta)$}

\subsection{Construction}

For any $n\in\N$, any elliptic curve $\E=\C/\Gamma$, and any point $\eta\in\E$ 
we construct a graded associative algebra 
$Q_n(\E,\eta)=\C\oplus F_1\oplus  F_2\oplus\dots$, where $F_1=\Theta_{n,c}(\Gamma)$ 
and $F_\alpha=S^\alpha\Theta_{n,c+(\alpha-1)n}(\Gamma)$. By construction, 
$\dim  F_\alpha=\frac{n(n+1)\dots(n+\alpha-1)}{\alpha!}$. It is clear that the 
space $F_\alpha$ can be realized as the space of holomorphic symmetric functions 
of $\alpha$ variables $\{f(z_1,\dots,z_\alpha)\}$ such that 
\begin{equation}
\begin{aligned}
f(z_1+1,z_2,\dots,z_\alpha)&=f(z_1,\dots,z_\alpha),\\
f(z_1+\tau,z_2,\dots,z_\alpha)&=(-1)^ne^{-2\pi i(nz_1-c-(\alpha-1)n)}f(z_1,\dots,z_\alpha).
\end{aligned}
\end{equation}

For $f\in  F_\alpha$ and $g\in  F_\beta$ we define the symmetric function $f*g$ 
of $\alpha+\beta$ variables by the formula 
\begin{multline*}
f*g(z_1,\dots,z_{\alpha+\beta})=\frac1{\alpha!\beta!}
\sum_{\sigma\in S_{\alpha+\beta}}
f(z_{\sigma_1},\dots,z_{\sigma_\alpha})
g(z_{\sigma_{\alpha+1}}-2\alpha\eta,\dots,z_{\sigma_{\alpha+\beta}}-
2\alpha\eta)\times\\
\times\prod_{\begin{subarray}{c}1\le                i\le\alpha\\\alpha+1\le
j\le\alpha+\beta\end{subarray}}
\frac{\theta(z_{\sigma_i}-z_{\sigma_j}-n\eta)}
{\theta(z_{\sigma_i}-z_{\sigma_j})}.
\end{multline*}
In particular, for $f,g\in F_1$ we have
$$
f*g(z_1,z_2)=f(z_1)g(z_2-2\eta)\frac{\theta(z_1-z_2-n\eta)}
{\theta(z_1-z_2)}+f(z_2)g(z_1-2\eta)\frac{\theta(z_2-z_1-n\eta)}
{\theta(z_2-z_1)}.
$$
Here $\theta(z)$ is a theta function of order one (see Appendix A).

\begin{prop}
If $f\in  F_\alpha$ and  $g\in F_\beta$, then $f*g\in F_{\alpha+\beta}$. The 
operation $*$ defines an associative multiplication on the space 
$\oplus_{\alpha\ge0}F_\alpha$
\end{prop}

\begin{proof}
Let us show that $f*g\in F_{\alpha+\beta}$. It immediately follows from the 
assumptions (9) concerning $f$ and $g$ and also from the properties of $\theta(z)$ 
(see Appendix A) that every summand in the formula for $f*g$ satisfies condition 
(9) for $F_{\alpha+\beta}$. Hence, $f*g$ is a meromorphic symmetric function 
satisfying condition (9). This function can have a pole of order not exceeding 
one on the diagonals $z_i-z_j=0$ and also for $z_i-z_j\in\Gamma$ because $\theta(z)$ 
has zeros for $z\in\Gamma$. However, the order of a pole of a symmetric function 
on the diagonal must be even. This implies that the function $f*g$ is holomorphic 
for $z_i=z_j$, and it follows from (9) that $f*g$ is holomorphic for $z_i-z_j\in\Gamma$ 
as well.

One can immediately see that the multiplication $*$ is associative.
\end{proof}

\subsection{Main properties of the algebra $Q_n(\E,\eta)$}

By construction, the dimensions of the graded components of the algebra   
$Q_n(\E,\eta)$ coincide with those for the polynomial ring in $n$ variables. For 
$\eta=0$ the formula for $f*g$ becomes
$$
f*g(z_1,\dots,z_{\alpha+1})=\frac1{\alpha!\beta!}\sum_{\sigma\in
S_{\alpha+\beta}}f(z_{\sigma_1},\dots,z_{\sigma_\alpha})
g(z_{\sigma_{\alpha+1}},\dots,z_{\sigma_{\alpha+\beta}}).
$$
This is the formula for the ordinary product in the algebra $S^*\Theta_{n,c}(\Gamma)$, 
that is, in the polynomial ring in $n$ variables. Therefore, for a fixed elliptic 
curve $\E$ (that is, for a fixed modular parameter $\tau$) the family of algebras 
$Q_n(\E,\eta)$ is a deformation of the polynomial ring. In particular (see 
Appendix D.2), there is a Poisson algebra, which we denote by $q_n(\E)$. One can 
readily obtain the formula for the Poisson bracket on the polynomial ring from 
the formula for $f*g$ by expanding the difference $f*g-g*f$ in the Taylor series 
with respect to $\eta$. It follows from the semicontinuity arguments that the 
algebra $Q_n(\E,\eta)$ with generic $\eta$ is determined by $n$ generators and 
$\frac{n(n-1)}2$ quadratic relations. One can prove (see \S2.6) that this is the 
case if $\eta$ is not a point of finite order on $\E$, that is, $N\eta\not\in\Gamma$ 
for any $N\in\N$.

The space $\Theta_{n,c}(\Gamma)$ of the generators of the algebra $Q_n(\E,\eta)$ 
is endowed with an action of a finite group $\wt{\Gamma_n}$ which is a central 
extension of the group $\Gamma/n\Gamma$ of points of order $n$ on the curve $\E$ 
(see Appendix A). It immediately follows from the formula for the product $*$ 
that the corresponding transformations of the space $F_\alpha=S^\alpha\Theta_{n,c}(\Gamma)$ 
are automorphisms of the algebra $Q_n(\E,\eta)$.

\subsection{Bosonization of the algebra $Q_n(\E,\eta)$}

The main approach to obtain representations of the algebra $Q_n(\E,\eta)$ is to 
construct homomorphisms from this algebra to other algebras with simple structure 
(close to Weil algebras) which have a natural set of representations. These 
homomorphisms are referred to as bosonizations, by analogy with the known 
constructions of quantum field theory.

Let $B_{p,n}(\eta)$ be a $\Z^p$-graded algebra whose space of degree 
$(\alpha_1,\dots,\alpha_p)$ is of the form  
$\{f(u_1,\dots,u_p)\e_1^{\alpha_1}\dots\e_p^{\alpha_p}\}$, where $f$ ranges over 
the meromorphic functions of $p$ variables and $\e_1,\dots,\e_p$ are elements of 
the algebra $B_{p,n}(\eta)$. Let $B_{p,n}(\eta)$ be generated by the space of 
meromorphic functions $f(u_1,\dots,u_p)$ and by the elements $\e_1,\dots,\e_p$ 
with the defining relations 
\begin{equation}
\begin{gathered}
\e_\alpha
f(u_1,\dots,u_p)=
f(u_1-2\eta,\dots,u_\alpha+(n-2)\eta,\dots,u_p-2\eta)\e_\alpha,\\
\e_\alpha\e_\beta=\e_\beta\e_\alpha,\quad
f(u_1,\dots,u_p)g(u_1,\dots,u_p)=g(u_1,\dots,u_p)f(u_1,\dots,u_p)
\end{gathered}
\end{equation}

We note that the subalgebra of $B_{p,n}(\eta)$ consisting of the elements of 
degree $(0,\dots,0)$ is the commutative algebra of all meromorphic functions of 
$p$ variables with the ordinary multiplication.

\begin{prop}
Let $\eta\in\E$ be a point of infinite order. For any $p\in\N$ there is a 
homomorphism $\phi_p\colon  Q_n(\E,\eta)\to B_{p,n}(\eta)$ that acts on the 
generators of the algebra $Q_n(\E,\eta)$ by the formula \emph:
\begin{equation}
\phi_p(f)=\sum_{1\le\alpha\le
p}\frac{f(u_\alpha)}{\theta(u_\alpha-u_1)\dots\theta(u_\alpha-u_p)}
\e_\alpha.
\end{equation}
Here $f\in\Theta_{n,c}(\Gamma)$ is a generator of $Q_n(\E,\eta)$ and the product 
in the denominator is of the form 
$\prod_{i\ne\alpha}\theta(u_\alpha-u_i)$.
\end{prop}

\begin{proof}
We write $\xi_\alpha=\frac1{\theta(u_\alpha-u_1)\dots\theta(u_\alpha-u_p)}
\e_\alpha$. It is clear that the elements $\xi_1,\dots,\xi_p$ together with the 
space of meromorphic functions $\{f(u_1,\dots,u_p)\}$ generate the algebra 
$B_{p,n}(\eta)$. The relations (10) become 
\begin{gather*}
\xi_\alpha
f(u_1,\dots,u_p)=
f(u_1-2\eta,\dots,u_\alpha+(n-2)\eta,\dots,u_p-2\eta)\xi_\alpha\\
\xi_\alpha\xi_\beta=-\frac{e^{2\pi
i(u_\beta-u_\alpha)}\theta(u_\alpha-u_\beta+n\eta)}
{\theta(u_\beta-u_\alpha+n\eta)}\xi_\beta\xi_\alpha
\end{gather*}
The formula (11) can be represented as 
\begin{equation}
\phi_p(f)=\sum_{1\le\alpha\le p}f(u_\alpha)\xi_\alpha.
\end{equation}
Using (12) and the formula for the multiplication in the algebra $Q_n(\E,\eta)$ 
and assuming that $\phi_p$ is a homomorphism, one can readily evaluate the 
extension of the map $\phi_p$ to the entire algebra. For instance, in the grading 
$2$ we have
\begin{align*}
\phi_p(f*g)&=\sum_{1\le\alpha\le p}f(u_\alpha)\xi_\alpha\cdot
\sum_{1\le\beta\le p}g(u_\beta)\xi_\beta=
\sum_{1\le\alpha,\beta\le                     p}f(u_\alpha)\xi_\alpha
f(u_\beta)\xi_\beta=\\
&=\sum_{\begin{subarray}{c}1\le\alpha,\beta\le
p\\\alpha\ne\beta\end{subarray}}
f(u_\alpha)g(u_\beta-2\eta)\xi_\alpha\xi_\beta+
\sum_{1\le\alpha\le p}f(u_\alpha)g(u_\alpha+(n-2)\eta)\xi_\alpha^2.
\end{align*}
The first sum is 
\begin{align*}
\sum_{1\le\alpha<\beta\le
p}&(f(u_\alpha)g(u_\beta-2\eta)\xi_\alpha\xi_\beta+
f(u_\beta)g(u_\alpha-2\eta)\xi_\beta\xi_\alpha)=
\end{align*}
\begin{align*}
=\sum_{1\le\alpha<\beta\le
p}\left(f(u_\alpha)g(u_\beta-2\eta)\xi_\alpha\xi_\beta-
f(u_\beta)g(u_\alpha-2\eta)\frac{e^{2\pi
i(u_\alpha-u_\beta)}\theta(u_\beta-u_\alpha-n\eta)}
{\theta(u_\alpha-u_\beta+n\eta)}\xi_\alpha\xi_\beta\right)=
\end{align*}
\begin{align*}
=\sum_{1\le\alpha<\beta\le
p}\frac{\theta(u_\alpha-u_\beta)}{\theta(u_\alpha-u_\beta-n\eta)}\times
\end{align*}
\begin{align*}
\times\left(f(u_\alpha)g(u_\beta-2\eta)
\frac{\theta(u_\alpha-u_\beta-n\eta)}{\theta(u_\alpha-u_\beta)}+
f(u_\beta)g(u_\alpha-2\eta)\frac{\theta(u_\beta-u_\alpha-n\eta)}
{\theta(u_\beta-u_\alpha)}\right)\xi_\alpha\xi_\beta=
\end{align*}
\begin{align*}
=\sum_{1\le\alpha<\beta\le
p}\frac{\theta(u_\alpha-u_\beta)}{\theta(u_\alpha-u_\beta-n\eta)}
f*g(u_\alpha,u_\beta)\xi_\alpha\xi_\beta\quad\text{ }
\end{align*}
Moreover, $f(u_\alpha)g(u_\alpha+(n-2)\eta)=
\frac{\theta(-n\eta)}{\theta(-2n\eta)}f*g(u_\alpha,u_\alpha+n\eta)$. We finally 
obtain
\medskip
\begin{equation}
\phi_p(f*g)=
\end{equation}
$$=\sum_{1\le\alpha<\beta\le
p}\frac{\theta(u_\alpha-u_\beta)}{\theta(u_\alpha-u_\beta-n\eta)}
f*g(u_\alpha,u_\beta)\xi_\alpha\xi_\beta+\frac{\theta(-n\eta)}
{\theta(-2n\eta)}\sum_{1\le\alpha\le
p}f*g(u_\alpha,u_\alpha+n\eta)\xi_\alpha^2$$
We see that the map $\phi_p$ can be extended to the quadratic part of the algebra 
$Q_n(\E,\eta)$ because the right-hand side of (13) depends on $f*g$ only but not 
on $f$ and $g$ separately. Thus implies the assertion for generic $\eta$ because 
in this case the algebra  $Q_n(\E,\eta)$ is defined by quadratic relations. To 
prove a more exact assertion (for the case in which $\eta$ is a point of infinite 
order), one must continue the above calculation. We obtain the following formula: 
if $f\in F_\alpha$, then 
\begin{equation}
\begin{aligned}
\phi_p(f)&=\sum_{\begin{subarray}{c}i_1,\dots,i_p\ge0,\\ i_1+\ldots+i_p=\alpha\end{subarray}}
B_{i_1,\dots,i_p}f(u_1,u_1+n\eta,\dots,u_1+(i_1-1)n\eta,\\
&\qquad\qquad\qquad\qquad u_2,u_2+n\eta,\dots,u_2+(i_2-1)n\eta,\dots)
\xi_1^{i_1}\dots\xi_p^{i_p},\\
&\text{where}\quad
B_{i_1,\dots,i_p}=\prod_{\begin{subarray}{c}1\le\lambda\le\lambda'\le
p\\0\le\mu<i_\lambda\\0\le\mu'<i_{\lambda'}\\\mu<\mu'\    \text{for}\
\lambda=\lambda'\end{subarray}}
\frac{\theta(u_\lambda+n\mu\eta-u_{\lambda'}-n\mu'\eta)}
{\theta(u_\lambda+n\mu\eta-u_{\lambda'}-n\mu'\eta-n\eta)}.
\end{aligned}
\end{equation}
This product can be represented as $\prod_{1\le
i<j<p}\frac{\theta(v_i-v_j)}{\theta(v_i-v_j-n\eta)}$, where 
$(v_1,\dots,v_p)=(u_1,u_1+n\eta,\dots)$ are the arguments of the function $f$ in 
the formula (14) for $\xi_1^{i_1}\dots\xi_p^{i_p}$.

The formula (14) makes sense if $\eta$ is a point of infinite order, and in 
this case the direct calculation shows that $\phi_p$ is a homomorphism.
\end{proof}

\subsection{Representations of the algebras $Q_n(\E,\eta)$}

The formula (14) shows that the image of the homomorphism $\phi_p$ is contained 
in the subalgebra $B_{p,n}^{\reg}(\eta)\subset  B_{p,n}(\eta)$ consisting of the 
elements $\sum_{\alpha_1,\dots,\alpha_p}f_{\alpha_1,\dots,\alpha_p}
\e^{\alpha_1}\dots\e_p^{\alpha_p}$, where the functions $f_{\alpha_1,\dots,\alpha_p}$ 
are holomorphic outside the divisors of the form 
$u_i-u_j-\lambda n\eta\in\Gamma$, $\lambda\in\Z$. 

Let $v_1,\dots,v_p\in\C$ be such that $v_i-v_j-\lambda
n\eta\not\in\Gamma$ for $\lambda\in\Z$. We construct a representation  
$M_{v_1,\dots,v_p}$ of the algebra $B_{p,n}^{\reg}(\eta)$ as follows. Let the 
representation $M_{v_1,\dots,v_p}$ have a basis 
$\{w_{\alpha_1,\dots,\alpha_p};\alpha_1,\dots\alpha_p\in\Z_{\ge0}\}$ in which the 
elements $\e_1,\dots,\e_p$ act by the rule $\e_iw_{\alpha_1,\dots,\alpha_p}=
w_{\alpha_1,\dots,\alpha_i+1,\dots,\alpha_p}$. Thus, 
$w_{\alpha_1,\dots,\alpha_p}=\e_1^{\alpha_1}\dots\e_p^{\alpha_p}w$, where 
$w=w_{0,\dots,0}$. The action of the commutative subalgebra of $B_{p,n}^{\reg}(\eta)$ 
consisting of the elements of degree $0$ is diagonal in this basis. We set 
$fw=f(v_1-(n-2)\eta,\dots,v_p-(n-2)\eta)w$, and hence $fw_{\alpha_1,\dots,\alpha_p}=
f(v_1+(2\alpha_1+\ldots+2\alpha_p-n\alpha_1-(n-2))\eta,\dots,
v_p+(2\alpha_1+\ldots+2\alpha_p-n\alpha_p-(n-2))\eta)
w_{\alpha_1,\dots,\alpha_p}$. It is clear that these formulas really define a 
representation of the algebra $B_{p,n}^{\reg}(\eta)$ in the space $M_{v_1,\dots,v_p}$, 
and, thanks to the homomorphism $\phi_p$, we have a representation of the algebra 
$Q_n(\E,\eta)$ as well. One can readily see that the space $M_{v_1,\dots,v_p}$ 
admits a basis $\{v_{\alpha_1,\dots,\alpha_p};\alpha_1,\dots,\alpha_p\in\Z_{\ge0}\}$, 
in which the action of the generators of the algebra $Q_n(\E,\eta)$ can be 
represented in the following form: if $f\in\Theta_n(\Gamma)$, then 
\begin{equation}
fv_{\alpha_1,\dots,\alpha_p}=\sum_{1\le                          i\le
p}\frac{f(v_i+(2\alpha_1+\ldots+2\alpha_p-n\alpha_i)\eta)}
{\theta(v_i-v_1-n(\alpha_i-\alpha_1)\eta)\dots
\theta(v_i-v_p-n(\alpha_i-\alpha_p)\eta)}
v_{\alpha_1,\dots,\alpha_i+1,\dots,\alpha_p}.
\end{equation}
The vectors $v_{\alpha_1,\dots,\alpha_p}$ are proportional to the vectors 
$w_{\alpha_1,\dots,\alpha_p}$. In particular, for $p=1$ we obtain modules $M_v$ 
with a basis $\{v_\alpha;\alpha\in\Z_{\ge0}\}$ and the action 
$fv_\alpha=f(v-(n-2)\alpha\eta)v_{\alpha+1}$. Thus, the algebra $Q_n(\E,\eta)$ 
has a family of linear modules parametrized by the elliptic curve $\E\subset\CP^{n-1}$, 
where the embedding is carried out by theta functions of order $n$.

\subsection{Symplectic leaves}

We recall that $Q_n(\E,0)$ is the polynomial ring $S^*\Theta_{n,c}(\Gamma)$. For 
a fixed elliptic curve $\E=\C/\Gamma$ we obtain the family of algebras $Q_n(\E,\eta)$, 
which is a flat deformation of the polynomial ring. We denote the corresponding 
Poisson algebra by $q_n(\E)$. We obtain a family of Poisson algebras, depending 
on $\E$, that is, on the modular parameter $\tau$. Let us study the symplectic 
leaves of this algebra. To this end, we note that, when passing to the limit as 
$\eta\to0$, the homomorphism $\phi_p$ of associative algebras gives a homomorphism 
of Poisson algebras. Namely, let us denote by $b_{p,n}$ the Poisson algebra formed 
by the elements 
$\sum_{\alpha_1,\dots,\alpha_p\ge0}
f_{\alpha_1,\dots,\alpha_p}(u_1,\dots,u_p)
\e_1^{\alpha_1}\dots\e_p^{\alpha_p}$, where $f_{\alpha_1,\dots,\alpha_p}$ are 
meromorphic functions and the Poisson bracket is 
$$
\{u_\alpha,u_\beta\}=\{\e_\alpha,\e_\beta\}=0;\quad
\{\e_\alpha,u_\beta\}=-2\e_\alpha;\quad
\{\e_\alpha,u_\alpha\}=(n-2)\e_\alpha,
$$
where $\alpha\ne\beta$.

The following assertion results from Proposition 6 in the limit as $\eta\to0$.

\begin{prop}
There is a Poisson algebra homomorphism $\psi_p\colon
q_n(\E)\to b_{p,n}$ given by the following formula: if $f\in\Theta_n(\Gamma)$, 
then
$\psi_p(f)=\sum_{1\le\alpha\le
p}\frac{f(u_\alpha)}{\theta(u_\alpha-u_1)\dots\theta(u_\alpha-u_p)}
\e_\alpha$.
\end{prop}

Let $\{\theta_i(u);i\in\Z/n\Z\}$ be a basis of the space 
$\Theta_{n,c}(\Gamma)$ and let $\{x_i;i\in\Z/n\Z\}$ be the corresponding basis 
in the space of elements of degree one in the algebra $Q_n(\E,\eta)$ (this space 
is isomorphic to $\Theta_{n,c}(\Gamma)$). For an elliptic curve $\E\subset\CP^{n-1}$ 
embedded by means of theta functions of order $n$ (this is the set of points with 
the coordinates $(\theta_0(z):\ldots:\theta_{n-1}(z))$) we denote by $C_p\E$ the 
variety of $p$-chords, that is, the union of projective spaces of dimension $p-1$ 
passing through $p$ points of $\E$. Let $K(C_p\E)$ be the corresponding homogeneous 
manifold in $\C^n$. It is clear that $K(C_p\E)$ consists of the points with the 
coordinates 
$x_i=\sum_{1\le\alpha\le
p}\frac{\theta_i(u_\alpha)}{\theta(u_\alpha-u_1)\dots
\theta(u_\alpha-u_p)}\e_\alpha$, where $u_\alpha,\e_\alpha\in\C$.

Let $2p<n$. Then one can show that $\dim K(C_p\E)=2p$ and $K(C_{p-1}\E)$ is the 
manifold of singularities of $K(C_p\E)$. It follows from Proposition 7 and from 
the fact that the Poisson bracket is non-degenerate on $b_{p,n}$ for $2p<n$ and 
$\e_\alpha\ne0$ that the non-singular part of the manifold $K(C_p\E)$ is a $2p$-
dimensional symplectic leaf of the Poisson algebra $q_n(\E)$.

Let $n$ be odd. One can show that the equation defining the manifold $K(C_{\frac{n-1}2}\E)$ 
is of the form $C=0$, where $C$ is a homogeneous polynomial of degree $n$ in the 
variables $x_i$. This polynomial is a central function of the algebra $q_n(\E)$. 

Let $n$ be even. The manifold $K(C_{\frac{n-2}2}\E)$ is defined by equations 
$C_1=0$ and $C_2=0$, where $\deg C_1=\deg C_2=n/2$. The polynomials $C_1$ and $C_2$ 
are central in the algebra $q_n(\E)$.

\subsection{Free modules, generations, and relations}

Let $\eta$ be a point of infinite order.

\begin{prop}
Let numbers $v_1,\dots,v_n\in\C$ be in general position. Then the module  
$M_{v_1,\dots,v_n}$ is generated by $v_{0,\dots,0}$ and is free over $Q_n(\E,\eta)$.
\end{prop}

\begin{proof}
By construction, the dimensions of graded components of $M_{v_1,\dots,v_n}$ 
coincide with those of the algebra $Q_n(\E,\eta)$. Let us show that the module is 
generated by the vector $v=v_{0,\dots,0}$. Let 
 
$$f_i=\prod_{\alpha\ne
i}\theta(z-v_\alpha)\cdot(\theta(z+v_1+\ldots+v_n-v_i-c).$$
It is clear that $f_i\in\Theta_{n,c}(\Gamma)$ for $1\le i\le n$. Therefore, the 
$f_i$s are elements of degree $1$ of the algebra $Q_n(\E,\eta)$. It follows from 
the formula (15) that $f_iv$ is non-zero and proportional to 
$v_{0,\dots,1,\dots,0}=\e_iv$. Similarly, one can readily construct elements 
$f_{i;\alpha_1,\dots,\alpha_n}\in\Theta_{n,c}(\Gamma)$ such that 
$f_{i;\alpha_1,\dots,\alpha_n}v_{\alpha_1,\dots,\alpha_n}$ is non-zero and 
proportional to $v_{\alpha_1,\dots,\alpha_i+1,\dots,\alpha_n}$. Namely,
$f_{i;\alpha_1,\dots,\alpha_n}=\prod_{\beta\ne
i}\theta(z-v_\beta-(2\alpha_1+\ldots+2\alpha_n-n\alpha_\beta)\eta)
\cdot\theta(z+v_1+\ldots+v_n-v_i+(n-2)(\alpha_1+\ldots+\alpha_n)\eta
-c)$.
Thus, all elements $v_{\alpha_1,\dots,\alpha_n}$ are obtained from $v$ by the action 
of elements of degree one in $Q_n(\E,\eta)$.
\end{proof}

\begin{prop}
The algebra $Q_n(\E,\eta)$ is presented by $n$ generators and $\frac{n(n-1)}2$ 
quadratic relations.
\end{prop}

\begin{proof}
It follows from the proof of Proposition 8 that the algebra $Q_n(\E,\eta)$ is generated 
by the elements of degree one. It is clear from the construction of the elements 
$f_{i;\alpha_1,\dots,\alpha_n}$ that these elements admit quadratic relations of 
the form
\begin{equation}
f_{j;\alpha_1,\dots,\alpha_i+1,\dots,\alpha_n}
f_{i;\alpha_1,\dots,\alpha_n}=
c_{i,j;\alpha_1,\dots,\alpha_n}
f_{i;\alpha_1,\dots,\alpha_j+1,\dots,\alpha_n}
f_{j;\alpha_1,\dots,\alpha_n},
\end{equation}
where $c_{i,j;\alpha_1,\dots,\alpha_n}\in\C^*$. To prove this relation, one must 
apply it to the vector $v_{\alpha_1,\dots,\alpha_n}$. Let us show that these 
quadratic relations imply the other ones. Let a relation be of the form 
$\sum_\alpha a_t^{(\alpha)}a_{t-1}^{(\alpha)}\dots a_1^{(\alpha)}=0$. We expand 
the element $a_1^{(\alpha)}$ in the basis $\{f_i\}$. The relation becomes 
$\sum_{\beta,i}b_t^{(\beta)}\dots b_2^{(\beta)}f_i=0$. Let us now expand $b_2^{(\beta)}$ 
in the basis $\{f_{i;0,\dots,1,\dots,0}\}$, where $1$ stands at the $i$th place. 
Continuing this procedure, we eventually represent the relation in the form $\sum
c_{i_1,\dots,i_t}f_{i_1;\alpha_1,\dots,\alpha_n}
f_{i_2;\alpha_1,\dots,\alpha_{i_2}-1,\dots,\alpha_n}\dots f_{i_t}=0$. It is clear 
that this relation follows from the relations (16).
\end{proof}

\begin{prop}
The relations in the algebra $Q_n(\E,\eta)$ are of the form 
\begin{equation}
\sum_{r\in\Z/n\Z}\frac1{\theta_{j-i-r}(-\eta)\theta_r(\eta)}x_{j-r}
x_{i+r}=0,\quad\text{сту}\ i\ne j;\ i,j\in\Z/n\Z.
\end{equation}
\end{prop}

\begin{proof}
Let us apply the formula for the multiplication in the algebra $Q_n(\E,\eta)$ 
(see \S2.1). Since $x_i=\theta_i(z)$, the relations (17) becomes 
\begin{multline*}
\sum_{r\in\Z/n\Z}\frac1{\theta_{j-i-r}(-\eta)\theta_r(\eta)}
\biggl(\theta_{j-r}(z_1)\theta_{i+r}(z_2-2\eta)
\frac{\theta(z_1-z_2-n\eta)}{\theta(z_1-z_2)}+\\
+\theta_{j-r}(z_2)\theta_{i+r}(z_1-2\eta)
\frac{\theta(z_2-z_1-n\eta)}{\theta(z_2-z_1)}\biggr)=0.
\end{multline*}
This relation immediately follows from the relation (30) (see Appendix A).
\end{proof}

\section{Main properties of the algebra $Q_{n,k}(\E,\eta)$}

We again assume that $\E=\C/\Gamma$ is an elliptic curve and $\eta\in\E$. Let 
$n$ and $k$ be coprime positive integers such that $1\le k<n$. Let us present 
the algebra $Q_{n,k}(\E,\eta)$ by the generators $\{x_i;i\in\Z/n\Z\}$ and the 
relations 
\begin{equation}
\sum_{r\in\Z/n\Z}\frac{\theta_{j-i+r(k-1)}(0)}
{\theta_{kr}(\eta)\theta_{j-i-r}(-\eta)}x_{j-r}x_{i+r}=0.
\end{equation}
As is known (see \S4), this is a PBW-algebra for generic $\E$ and $\eta$. We 
conjecture that this holds for any $\E$ and $\eta$. For generic $\E$ and $\eta$ 
the centre of the algebra $Q_{n,k}(\E,\eta)$ is the polynomial ring in $c=\NOD(n,k+1)$ 
elements of degree $n/c$ (see [20]). Hypothetically, this is the case for any 
$\E$ and $\eta$, where $\eta$ is a point of infinite order. If $\eta\in\E$ is a 
point of finite order, then the algebra $Q_{n,k}(\E,\eta)$ is finite-dimensional 
over its centre (see [37]). The following properties can readily be verified:

\smallskip

1) $Q_{n,k}(\E,0)=\C[x_1,\dots,x_n]$ is commutative;

2) $Q_{n,n-1}(\E,\eta)=\C[x_1,\dots,x_n]$ is commutative for any $\eta$;

3) $Q_{n,k}(\E,\eta)\simeq    Q_{n,k'}(\E,\eta)$, where $kk'\equiv1\pmod n$;

4) the maps $x_i\mapsto x_{i+1}$ and $x_i\mapsto\epsilon^ix_i$ (where $\epsilon$ 
is a primitive root of unity of degree $n$) define automorphisms of the algebra 
$Q_{n,k}(\E,\eta)$.

\smallskip

It follows from the results of \S5 (see Proposition 11) that the space of generators 
of the algebra $Q_{n,k}(\E,\eta)$ is naturally isomorphic to the space of theta 
functions $\Theta_{n/k}(\Gamma)$ (see Appendix B). Moreover, this space of generators 
is dual to the space $\Theta_{n/n-k}(\Gamma)$  (see Proposition 14). For a 
description of this duality between the spaces of theta functions, see Appendix C.

The algebra $Q_{n,k}(\E,\eta)$ is not a Hopf algebra and admits no comultiplications. 
However, there are homomorphisms of the algebra $Q_{n,k}(\E,\eta)$ to tensor 
products of other algebras of this kind (see [35; \S3]). To describe these 
homomorphisms, we need the notation of Appendix B. Moreover, we denote by 
$L_m(\E,\eta)=\C\oplus\Theta_{m,0}(\Gamma)\oplus
\Theta_{2m,m\eta}(\Gamma)\oplus\dots$ the $\Z_{\ge0}$-graded algebra with the 
multiplication $*$ given by the formula $f*g(z)=f(z+\beta\eta)g(z)$, where $\beta$ 
is a power of $g$. As we know from \S2, $L_n(\E,(n-2)\eta)$ is a quotient algebra 
of $Q_n(\E,\eta)$.

Let $A$ be an associative algebra and let $G\subset\Aut  A$. We denote by 
$A^G\subset  A$ the subalgebra consisting of the elements invariant with respect to 
$G$.

There are the following algebra homomorphisms:

\smallskip

a)
$Q_{n,k}(\E,\eta)\to
\left(L_{kn}\left(\E,\frac{n-k-1}k\eta\right)\otimes
Q_{k,\l}\left(\E,\frac   nk\eta\right)\right)^{\wt{\Gamma_k}}$,
where $\l=d(n_3,\dots,n_p)$ and the generators are taken to elements of bidegree 
$(1,1)$.

b)
$Q_{n,k}(\E,\eta)\to
\left(L_{nk'}\left(\E,\frac{n-k'-1}{k'}\eta\right)\otimes
Q_{k',\l'}\left(\E,\frac
n{k'}\eta\right)\right)^{\wt{\Gamma_{k'}}}$,
where $\l'=d(n_1,\dots,n_{p-2})$ and the generators are taken to elements of bidegree 
$(1,1)$.

c)
$Q_{n,k}(\E,\eta)\to
\left(Q_{a,\alpha}\left(\E,\frac n\alpha\eta\right)\otimes
L_{abn}\left(\E,\frac{n-a-b}{ab}\eta\right)\otimes
Q_{b,\beta}\left(\E,\frac
nb\right)\right)^{\wt{\Gamma_{ab}}}$,
where $a=d(n_1,\dots,n_{i-1})$, $b=d(n_{i+1},\dots,n_p)$, $\alpha=d(n_1,\dots,n_{i-2})$, and  
$\beta=d(n_{i+2},\dots,n_p)$ for some $i$; the generators are taken to elements 
of multidegree $(1,1,1)$.

\smallskip

Let us describe the map c) geometrically (the description of the maps a) and b) 
is the same). Let $f(z_1,\dots,z_p)\in\Theta_{n/k}(\Gamma)$. For some $i$  ($1<i<p$) 
we choose a $z_i$, then $f$ (regarded as a function of $z_1,\dots,z_{i-1}$) 
belongs to a space isomorphic to $\Theta_{a/\alpha}(\Gamma)$. Similarly, when 
regarded as a function of $z_{i+1},\dots,z_p$, $f$ belongs to a space isomorphic to 
$\Theta_{b/\beta}(\Gamma)$. Thus, for a fixed $z_i$ we have 
$f\in\Theta_{a/\alpha}(\Gamma)\otimes\Theta_{b/\beta}(\Gamma)$. A family of linear 
maps $\Theta_{n/k}(\Gamma)\to\Theta_{a/\alpha}(\Gamma)\otimes
\Theta_{b/\beta}(\Gamma)$ arises. With regard to the dependence on $z_i$, we 
obtain a linear map $\Theta_{n/k}(\Gamma)\to\Theta_{a/\alpha}(\Gamma)\otimes
\Theta_{nab,0}(\Gamma)\otimes\Theta_{b/\beta}(\Gamma)$. The homomorphism c) corresponds 
to this map (the space of generators of the algebra $Q_{n,k}(\E,\eta)$ is 
$\Theta_{n/k}(\Gamma)$, the space of generators of the algebra 
$L_{nab}\left(\E,\frac{n-a-b}{ab}\eta\right)$ is $\Theta_{nab,0}(\Gamma)$, etc.).

\section{Belavin elliptic $R$-matrix and the algebra $Q_{n,k}(\E,\eta)$}

Let $V$ be a vector space of dimension $n$. For each $u\in\C$ we denote by $V(u)$ 
a vector space canonically isomorphic to $V$. Let $R$ be a meromorphic function 
of two variables with values in $\End(V\otimes V)$. It is convenient to regard 
$R(u,v)$ as a linear map 
$$
R(u,v)\colon V(u)\otimes V(v)\to V(v)\otimes V(u).
$$

We recall that by the Yang-Baxter equation one means the condition that the 
following diagram is commutative:

{\footnotesize

$$
\xymatrix{&V(v)\otimes V(u)\otimes V(w)\ar[r]^{1\otimes R(u,w)}&
V(v)\otimes V(w)\otimes V(u)\ar[dr]^{R(v,w)\otimes1}\\
V(u)\otimes V(v)\otimes V(w)\ar[ur]^{R(u,v)\otimes1}\ar[dr]^{1\otimes
R(v,w)}&&&
V(w)\otimes V(v)\otimes V(u)\\
&V(u)\otimes V(w)\otimes V(v)\ar[r]^{R(u,w)\otimes1}&
V(w)\otimes V(u)\otimes V(v)\ar[ur]_{1\otimes R(u,v)}}
$$}

A solution of the Yang-Baxter equation is called $R$-matrix.

Let $\{x_i;i=1,\dots,n\}$ be a basis in the space $V$ and let $\{x_i(u)\}$ be the 
corresponding basis in the space $V(u)$. Let $R_{ij}^{\alpha\beta}(u,v)$ be an 
entry of an $R$-matrix $R(u,v)$, that is, $R(u,v)\colon     x_i(u)\otimes     x_j(v)\to
R_{ij}^{\alpha\beta}(u,v)x_\beta(v)\otimes x_\alpha(u)$.

Let an $R$-matrix $R(u,v)$ satisfy the relation $R(u,v)R(v,u)=1$. By the 
Zamolodchikov algebra $Z_R$ one means the algebra with the generators 
$\{x_i(u);i=1,\dots,n;u\in\C\}$ and the defining relations 
\begin{equation}
x_i(u)x_j(v)=R_{ij}^{\alpha\beta}(u,v)x_\beta(v)x_\alpha(u).
\end{equation}
It is clear that the elements $\{x_{i_1}(u_1)\dots
x_{i_m}(u_m);1\le i_1,\dots,i_m\le n\}$ of the Zamolodchikov algebra are linearly 
independent for generic $u_1,\dots,u_m$. Thus, Zamolodchikov algebras are 
infinite-dimensional analogues of PBW-algebras. We recall that by the classical 
$r$-matrix one means a Poisson structure of the form 
$\{x_i(u),x_j(v)\}=r_{ij}^{\alpha\beta}(u,v)x_\alpha(u)x_\beta(v)$. It is clear 
that the Zamolodchikov algebra is a quantization of this Poisson structure if the 
$R$-matrix depends on an additional parameter $\hbar$ and the relations (19) are 
of the form 
$$
x_i(u)x_j(v)=x_j(v)x_i(u)+\hbar
r_{ij}^{\alpha\beta}(u,v)x_\beta(v)x_\alpha(u)+o(\hbar).
$$

The Yang-Baxter equation has elliptic solutions, which are referred to as Belavin 
$R$-matrices. Let $n,k\in\N$ be coprime and $1\le k<n$. For any $n$ and $k$ 
there is a two-parameter family of $R$-matrices $R_{n,k}(\E,\eta)$ depending on 
an elliptic curve $\E=\C/\Gamma$ and a point $\eta\in\E$. Namely,
$$
R_{n,k}(\E,\eta)(u-v)                       (x_i(u)\otimes
x_j(v))=\frac1{p(u-v)}\sum_{r\in\Z/n\Z}
\frac{\theta_{j-i+r(k-1)}(v-u+\eta)}
{\theta_{kr}(\eta)\theta_{j-i-r}(v-u)}x_{j-r}(v)\otimes x_{i+r}(u),
$$
where $p(u-v)=
\frac{\theta_1(0)\dots\theta_{n-1}(0)
\theta_0(v-u+\eta)\dots\theta_{n-1}(v-u+\eta)}
{\theta_0(\eta)\dots\theta_{n-1}(\eta)
\theta_0(v-u)\dots\theta_{n-1}(v-u)}$;     $i,j\in\Z/n\Z$.
One can readily see that $\det
R_{n,k}(\E,\eta)(u-v)=
\left(\frac{\theta_0(v-u-\eta)\dots\theta_{n-1}(v-u-\eta)}
{\theta_0(v-u+\eta)\dots\theta_{n-1}(v-u+\eta)}\right)
^{\frac{n(n-1)}2}$. 
Thus, the operator $R_{n,k}(\E,\eta)(-\eta)$ has a kernel. Let $L_{n,k}(\E,\eta)\subset
V\otimes  V$ and $L_{n,k}(\E,\eta)=\ker R_{n,k}(\E,\eta)(-\eta)$. According to [10] 
we have $\dim(L_{n,k}(\E,\eta))=\frac{n(n-1)}2$, and $L_{n,k}(\E,0)=\Lambda^2V$ 
for $\eta=0$. Let $Q_{n,k}(\E,\eta)=T^*V/(L_{n,k}(\E,\eta))$ be the algebra with 
the generators $\{x_i;i\in\Z/n\Z\}$ and the defining relations $L_{n,k}(\E,\eta)$. 
The dimensions of the graded components of the algebra $Q_{n,k}(\E,\eta)$ coincide 
with those of the polynomial ring $S^*V$ (see [10]). It follows from the formula 
for $R_{n,k}(\E,\eta)(u-v)$ that the defining relations in the algebra $Q_{n,k}(\E,\eta)$ 
are 
\begin{equation}
\sum_{r\in\Z/n\Z}\frac{\theta_{j-i+r(k-1)}(0)}
{\theta_{kr}(\eta)\theta_{j-i-r}(-\eta)}x_{j-r}x_{i+r}=0.
\end{equation}
In particular, we see that $Q_{n,1}(\E,\eta)=Q_n(\E,\eta)$ if $\eta\in\E$ is a 
point of infinite order.

It follows from the relations (20) that $Q_{n,k}(\E,\eta)\simeq
Q_{n,k'}(\E,\eta)$, where $kk'\equiv1\pmod  n$. Moreover, $Q_{n,n-1}(\E,\eta)$ is 
commutative for any $\E,\eta$. It is also clear that $Q_{n,k}(\E,0)$ is 
commutative.

\section{Algebras $Q_{n,k}(\E,\eta)$ and the exchange algebras}

\subsection{Homomorphisms of algebras $Q_{n,k}(\E,\eta)$ into dynamical exchange 
algebras}

The algebras $Q_{n,k}(\E,\eta)$ for arbitrary $k$ have representations similar to 
the homomorphisms $\phi_p$ related to the case $k=1$ (see \ \S2.3,  (11)). 
However, the structure of the algebra similar to $B_{p,n}(\eta)$ for $k=1$ turns 
out to be more complicated for $k>1$. Let $n/k=n_1-\frac1{n_2-\ldots-\frac1{n_q}}$ 
be the expansion of the number $n/k$ in the continued fraction in which $n_\alpha\ge2$ 
for $1\le\alpha\le q$. It is clear that such an expansion exists and is unique. 
We recall that $1\le    k<n$, where $n$ and $k$ are coprime. Let 
$d(m_1,\dots,m_t)=\det M$, where $M=(m_{ij})$ be a ($t\times t$) matrix with the 
entries $m_{ii}=m_i$,  $m_{i,i+1}=m_{i+1,i}=-1$, and $m_{ij}=0$ for $|i-j|>1$. 
For $t=0$ we set $d(\varnothing)=1$. It is clear that 
$n=d(n_1,\dots,n_q)$ and $k=d(n_2,\dots,n_q)$.

\begin{prop}
There is an algebra homomorphism of $Q_{n,k}(\E,\eta)$ into the algebra $C_{n,k}(\eta)$ 
generated by the commutative subalgebra $\{f(y_1,\dots,y_q)\}$ of holomorphic 
functions of $q$ variables and by an element $\e$ with the defining relations of 
the form $\e f(y_1,\dots,y_q)=f(y_1+\eta_1,\dots,y_q+\eta_q)\e$, where  
$\eta_\alpha=(d(n_1,\dots,n_q)-d(n_1,\dots,n_{\alpha-1})-
d(n_{\alpha+1},\dots,n_q))\eta$.  Moreover, 
$x_i\to w_i(y_1,\dots,y_q)\e$, where
$w_i\in\Theta_{n/k}(\Gamma)$ \emph(see Appendix B).
\end{prop}

This is a special case of Proposition 13 below.

The algebra $C_{n,k}(\eta)$ has a family of modules $L_{v_1,\dots,v_q}$ with a 
basis $\{v_\alpha;\alpha\in\Z_{\ge0}\}$ and with the action given by $\e
v_\alpha=v_{\alpha+1}$ and
$f(y_1,\dots,y_q)v_\alpha=
f(v_1-\alpha\eta_1,\dots,v_q-\alpha\eta_q)v_\alpha$. This implies the following 
assertion.

\begin{prop}
The algebra $Q_{n,k}(\E,\eta)$ hes a family of modules $L_{u_1,\dots,u_q}$ with 
a basis $\{v_\alpha;\alpha\in  Z_{\ge0}\}$ and with the action \emph:
$x_iv_\alpha=w_i(u_1-\alpha\eta_1,\dots,u_q-\alpha\eta_q)
v_{\alpha+1}$, where $w_i\in\Theta_{n/k}(\Gamma)$      and
$\eta_j=(d(n_1,\dots,n_q)-d(n_1,\dots,n_{j-1})-
d(n_{j+1},\dots,n_q))\eta$.
\end{prop}

In particular, we see that the algebra $Q_{n,k}(\E,\eta)$ has a family of linear 
modules depending on the point of $\E^q$, and the space of generators of the 
algebra $Q_{n,k}(\E,\eta)$ is isomorphic to $\Theta_{n/k}(\Gamma)$.

Let $C_{m_1,\dots,m_q;n,k}(\E,\eta)$ be the algebra generated by the commutative 
subalgebra $\{\phi(y_{1,1},\dots,y_{m_1,1};\dots;y_{1,q},\dots,y_{n_q,q})\}$, 
where $\phi$ are the meromorphic functions with the ordinary multiplication, and 
by elements $\{\e_{\alpha_1,\dots,\alpha_q};1\le\alpha_t\le m_t\}$. The defining 
relations for the algebra $C_{m_1,\dots,m_q;n,k}(\E,\eta)$ look as follows:
\begin{align*}
\e_{\alpha_1,\dots,\alpha_q}y_{\beta,\nu}&=
(y_{\beta,\nu}-(d(n_1,\dots,n_{\nu-1})+d(n_{\nu+1},\dots,n_q))\eta)
\e_{\alpha_1,\dots,\alpha_q},\quad\ \alpha_\nu\ne\beta,\\
\e_{\alpha_1,\dots,\alpha_q}y_{\alpha_\nu,\nu}&=
(y_{\alpha_\nu,\nu}+(n-d(n_1,\dots,n_{\nu-1})-d(n_{\nu+1},\dots,n_q))
\eta)\e_{\alpha_1,\dots,\alpha_q}.
\end{align*}
These relations mean that the $y_{\beta,\nu}$s are dynamical variables. This 
immediately implies the relations between the elements $\e_{\alpha_1,\dots,\alpha_q}$ 
and the meromorphic functions of the variables $y_{\beta,\nu}$. The remaining 
relations are quadratic in $\e_{\alpha_1,\dots,\alpha_q}$ with the coefficients 
depending on the dynamical variables $y_{\beta,\nu}$. The relations in "general 
position" are as follows:
\begin{equation}
\e_{\alpha_1,\dots,\alpha_q}\e_{\beta_1,\dots,\beta_q}=
\Lambda\e_{\beta_1,\dots,\beta_q}\e_{\alpha_1,\dots,\alpha_q}+
\sum_{1\le t\le q-1}\Lambda_{t,t+1}
\e_{\beta_1,\dots,\beta_t,\alpha_{t+1},\dots,\alpha_q}
\e_{\alpha_1,\dots,\alpha_t,\beta_{t+1},\dots,\beta_q},
\end{equation}
where $\alpha_1\ne\beta_1$, \dots, $\alpha_q\ne\beta_q$ and
$$\begin{aligned}
\Lambda&=\frac{e^{-2\pi in\eta}\theta(y_{\beta_1,1}-y_{\alpha_1,1})
\theta(y_{\beta_q,q}-y_{\alpha_q,q}+n\eta)}
{\theta(y_{\beta_1,1}-y_{\alpha_1,1}-n\eta)
\theta(y_{\beta_q,q}-y_{\alpha_q,q})},\\
\Lambda_{t,t+1}&=\frac{e^{-2\pi in\eta}\theta(n\eta)
\theta(y_{\beta_1,1}-y_{\alpha_1,1})}
{\theta(y_{\beta_1,1}-y_{\alpha_1,1}-n\eta)}\cdot
\frac{\theta(y_{\beta_t,t}+y_{\beta_{t+1},t+1}-
y_{\alpha_t,t}-y_{\alpha_{t+1},t+1})}
{\theta(y_{\beta_t,t}-y_{\alpha_t,t})
\theta(y_{\beta_{t+1},t+1}-y_{\alpha_{t+1},t+1})}.
\end{aligned}$$
The relations of non-general position occur if some $\alpha_\nu$s are equal to  
$\beta_\nu$s. These relations exist for any subset of the form 
$\{\psi+1,\dots,\psi+\phi\}$, where $0\le\psi$, $\psi+\phi\le q$, and 
$\alpha_\psi=\beta_\psi$               (if               $0<\psi$),
$\alpha_{\psi+\phi+1}=\beta_{\psi+\phi+1}$  (if $\psi+\phi<q$), and
$\alpha_{\psi+1}\ne\beta_{\psi+1}$,                            \dots,
$\alpha_{\psi+\phi}\ne\beta_{\psi+\phi}$. These relations are of the form
\begin{multline}
\e_{\mu_1,\dots,\mu_\psi,\alpha_1,\dots,\alpha_\phi,
\gamma_1,\dots,\gamma_p}
\e_{\mu_1',\dots,\mu_{\psi-1}',\mu_\psi,\beta_1,\dots,\beta_\phi,
\gamma_1,\gamma_2',\dots,\gamma_p'}=\\
=\Lambda\e_{\mu_1,\dots,\mu_\psi,\beta_1,\dots,\beta_\phi,
\gamma_1,\dots,\gamma_p}
\e_{\mu_1',\dots,\mu_\psi,\alpha_1,\dots,\alpha_\phi,
\gamma_1,\dots,\gamma_p'}+\\
+\sum_{1\le t<\phi}\Lambda_{t,t+1}\e_{\mu_1,\dots,\mu_\psi,
\beta_1,\dots,\beta_t,\alpha_{t+1},\dots,\alpha_\phi,
\gamma_1,\dots,\gamma_p}
\e_{\mu_1',\dots,\mu_\psi,\alpha_1,\dots,\alpha_t,
\beta_{t+1},\dots,\beta_\phi,\gamma_1,\dots,\gamma_p'}.
\end{multline}
Here $\alpha_1\ne\beta_1$,   \dots,    $\alpha_\phi\ne\beta_\phi$, and
$\phi+\psi+p=q$. The coefficients $\Lambda,\Lambda_{t,t+1}$ are defined by (21). 
If $\psi=p=0$, then this relation coincides with (21), that is, becomes a relation 
in general position.

We note that, if $\eta=0$, then the algebra $C_{m_1,\dots,m_q;n,k}(\E,0)$ is 
commutative and does not depend on the elliptic curve $\E$. For $q=1$ and $q=2$ 
the algebra $C_{m_1,\dots,m_q;n,k}(\E,0)$ is the polynomial ring in the variables 
$\{\e_{\alpha_1,\dots,\alpha_q}\}$ over the field of meromorphic functions of the 
variables $\{y_{\alpha,\beta}\}$. For $q>2$ the algebra $C_{m_1,\dots,m_q;n,k}(\E,0)$ 
is no longer a polynomial ring because additional relations occur. Namely, the 
relations (22) for $\eta=0$ become 
\begin{multline}
\e_{\mu_1,\dots,\mu_\psi,\alpha_1,\dots,\alpha_\phi,
\gamma_1,\dots,\gamma_p}
\e_{\mu_1',\dots,\mu_{\psi-1}',\mu_\psi,\beta_1,\dots,\beta_\phi,
\gamma_1,\gamma_2',\dots,\gamma_p'}=\\
=\e_{\mu_1,\dots,\mu_\psi,\beta_1,\dots,\beta_\phi,
\gamma_1,\dots,\gamma_p}
\e_{\mu_1',\dots,\mu_{\psi-1}',\mu_\psi,\alpha_1,\dots,\alpha_\phi,
\gamma_1,\gamma_2',\dots,\gamma_p'}.
\end{multline}
The algebra $C_{m_1,\dots,m_q;n,k}(\E,\eta)$ is a flat deformation of the function 
algebra on the manifold defined by the relations (23). Moreover, $y_{\alpha,\beta}$ 
are dynamical variables. It is clear that the manifold defined by the equations 
(23) is rational, and the general solution of these equations is 
$\e_{\alpha_1,\dots,\alpha_q}=\e_{\alpha_1,\alpha_2}^{(1)}
\e_{\alpha_2,\alpha_3}^{(2)}\dots\e_{\alpha_{q-1},\alpha_q}^{(q-1)}$,
where $\{\e_{\alpha,\beta}^{(t)}\}$ are independent variables. 

\begin{prop}
There is an algebra homomorphism 
$\phi\colon  Q_{n,k}(\E,\eta)\to
C_{m_1,\dots,m_q;n,k}(\E,\eta)$ that acts on the generators of the algebra 
$Q_{n,k}(\E,\eta)$ as follows: \emph:
\begin{equation}
\phi(x_i)=\sum_{\begin{subarray}{c}1\le\alpha_1\le
m_1\\\dots\\1\le\alpha_q\le m_q\end{subarray}}
w_i(y_{\alpha_1,1},\dots,y_{\alpha_q,q})
\e_{\alpha_1,\dots,\alpha_q}.
\end{equation}
Here $w_i\in\Theta_{n/k}(\Gamma)$.
\end{prop}

\begin{proof}
The algebra $Q_{n,k}(\E,\eta)$ is defined by the relations (18). Let us show that 
the images $\phi(x_i)$ satisfy the same relations. We have
\begin{multline*}
\sum_{r\in\Z/n\Z}\frac{\theta_{j-i+r(k-1)}(0)}
{\theta_{kr}(\eta)\theta_{j-i-r}(-\eta)}\phi(x_{j-r})\phi(x_{i+r})=\\
=\sum_{\begin{subarray}{c}r\in\Z/n\Z\\\\1\le\alpha_1,\beta_1\le
m_1\\\dots\\1\le\alpha_q,\beta_q\le m_q\end{subarray}}
\frac{\theta_{j-i+r(k-1)}(0)}{\theta_{kr}(\eta)\theta_{j-i-r}(-\eta)}
w_{j-r}(y_{\alpha_1,1},\dots,y_{\alpha_q,q})
\e_{\alpha_1,\dots,\alpha_q}
w_{i+r}(y_{\beta_1,1},\dots,y_{\beta_q,q})
\e_{\beta_1,\dots,\beta_q}.
\end{multline*}
Using the relations (21) and (22), we obtain an expression of the form 
$$
\sum_{\begin{subarray}{l}\alpha_1\le\beta_1;\\\alpha_2\le\beta_2\
\text{for}\
\alpha_1=\beta_1;\\
\dots\end{subarray}}
\psi_{\alpha_1,\dots,\alpha_q,\beta_1,\dots,\beta_q}
\e_{\alpha_1,\dots,\alpha_q}\e_{\beta_1,\dots,\beta_q}.
$$
We must prove that the coefficients are equal to $0$. Let us restrict ourselves 
to the case $\alpha_1\ne\beta_1$, \dots, $\alpha_q\ne\beta_q$. In this case we 
have by direct calculation
\begin{multline*}
\psi_{\alpha_1,\dots,\beta_q}=\sum_{r\in\Z/n\Z}
\frac{\theta_{j-i+r(k-1)}(0)}{\theta_{kr}(\eta)\theta_{j-i-r}(-\eta)}
\biggl(w_{j-r}(y_{\alpha_1,1},\dots,y_{\alpha_q,q})
w_{i+r}(y_{\beta_1,1}+t_1,\dots,y_{\beta_q,q}+t_q)+\\
+\frac{e^{-2\pi in\eta}\theta(y_{\alpha_1,1}-y_{\beta_1,1})
\theta(y_{\alpha_q,q}-y_{\beta_q,q}+n\eta)}
{\theta(y_{\alpha_1,1}-y_{\beta_1,1}-n\eta)
\theta(y_{\alpha_q,q}-y_{\beta_q,q})}
w_{j-r}(y_{\beta_1,1},\dots,y_{\beta_q,q})\times\\
\times w_{i+r}(y_{\alpha_1,1}+t_1,\dots,y_{\alpha_q,q}+t_q)+\\
+\sum_{1\le t\le q-1}\frac{e^{-2\pi in\eta}\theta(n\eta)
\theta(y_{\alpha_1,1}-y_{\beta_1,1})}
{\theta_{\alpha_1,1}-y_{\beta_1,1}-n\eta)}\cdot
\frac{\theta(y_{\alpha_t,t}+y_{\beta_{t+1},t+1}-y_{\beta_t,t}-
y_{\alpha_{t+1},t+1})}{\theta(y_{\alpha_t,t}-y_{\beta_t,t})
\theta(y_{\beta_{t+1},t+1}-y_{\alpha_{t+1},t+1})}\times\\
\times w_{j-r}(y_{\beta_1,1},\dots,y_{\beta_t,t},
y_{\alpha_{t+1},t+1},\dots,,y_{\alpha_q,q})\times\\
\times w_{i+r}(y_{\alpha_1,1}+t_1,\dots,y_{\alpha_t}+t_t,
y_{\beta_{t+1}}+t_{t+1},\dots,y_{\beta_q,q}+t_q)\biggr).
\end{multline*}
Here
$t_\alpha=-(d(n_1,\dots,n_{\alpha-1})+d(n_{\alpha+1},\dots,n_q))\eta$.

The equality $\psi_{\alpha_1,\dots,\beta_q}=0$ immediately follows from the 
identity (35) proved in Appendix B.
\end{proof}

\subsection{Homomorphism of the exchange algebra into the algebra
 $Q_{n,k}(\E,\eta)$}

Let $q'\in\N$ and $\mu_1,\dots,\mu_{q'},\mu\in\C$. We define an associative algebra 
$Y_{q'}(\E,\mu;\mu_1,\dots,\mu_{q'})$ as follows. The algebra 
$Y_{q'}(\E,\mu;\mu_1,\dots,\mu_{q'})$ is presented by the generators 
$\{\e(u_1,\dots,u_{q'});u_1,\dots,u_{q'}\in\C\}$ and the defining relations 
\begin{multline*}
\frac{\theta(v_1-u_1+\mu)}{\theta(v_1-u_1)} \e(u_1,\dots,u_{q'})
\e(v_1+\mu_1,\dots,v_{q'}+\mu_{q'})=\\
=\sum_{1\le t<q'}\frac{\theta(\mu)\theta(v_t-u_t+u_{t+1}-v_{t+1})}
{\theta(v_t-u_t)\theta(u_{t+1}-v_{t+1})}\times\\
\times\e(v_1,\dots,v_t,u_{t+1},\dots,u_{q'})
\e(u_1+\mu_1,\dots,u_t+\mu_t,v_{t+1}+\mu_{t+1},\dots,v_{q'}+\mu_{q'})
+\\
+\frac{\theta(v_{q'}-u_{q'}+\mu)}{\theta(v_{q'}-u_{q'})}
\e(v_1,\dots,v_{q'})\e(u_1+\mu_1,\dots,u_{q'}+\mu_{q'}).
\end{multline*}

For $\mu=\mu_1=\ldots=\mu_{q'}=0$ the algebra $Y_{q'}(\E,0;0,\dots,0)$ is the 
polynomial ring in infinitely many variables 
$\{\e(u_1,\dots,u_{q'});u_1,\dots,u_{q'}\in\C\}$. One can show that the algebra 
$Y_{q'}(\E,\mu;\mu_1,\dots,\mu_{q'})$ is a flat deformation of this polynomial 
ring. 

Let $\frac     n{n-k}=n_1'-\frac1{n_2'-\ldots-\frac1{n_{q'}'}}$ be an expansion 
in the continued fraction in which $n_\alpha'\ge2$ for $1\le\alpha\le q'$. It is 
clear that such an expansion exists and unique. For the relationship between the 
expansions in continued fractions of the numbers $\frac  nk$ and $\frac n{n-k}$, 
see Appendix C. 

\begin{prop}
There is an algebra homomorphism 
$$\psi\colon
Y_{q'}(\E,\mu;\mu_1,\dots,\mu_{q'})\to Q_{n,k}(\E,\eta),$$
 where $\mu=n\eta$ and
$\mu_\alpha=(d(n_1',\dots,n_{\alpha-1}')-d(n_{\alpha+1}',
\dots,n_{q'}')\eta$. This homomorphism is of the form \emph:
\begin{equation}
\psi\colon\e(u_1,\dots,u_{q'})\to\sum_{\alpha\in\Z/n\Z}
w_\alpha(u_1,\dots,u_{q'})x_{1-\alpha},
\end{equation}
where $\{x_i;i\in\Z/n\Z\}$ are the generators of the algebra $Q_{n,k}(\E,\eta)$, and 
$\{w_\alpha;\alpha\in\Z/n\Z\}$ is a basis in the space of theta functions 
$\Theta_{n/n-k}(\Gamma)$ (see Appendix B).
\end{prop}

\begin{proof}
Let us apply the map $\psi$ to the difference between the left- and right-hand 
sides of the relations in the algebra $Y_{q'}(\E,\mu;\mu_1,\dots,\mu_{q'})$. We 
must verify the resulting relation in the algebra $Q_{n,k}(\E,\eta)$. We have 
\begin{multline*}
\frac{\theta(v_1-u_1+\mu)}{\theta(v_1-u_1)}
\psi(\e(u_1,\dots,u_{q'}))\psi(\e(v_1+\mu_1,\dots,v_{q'}+\mu_{q'}))-
\\
-\sum_{1\le t<q'}\frac{\theta(\mu)\theta(v_t-u_t+u_{t+1}-v_{t+1})}
{\theta(v_t-u_t)\theta(u_{t+1}-v_{t+1})}\times\\
\times\psi(\e(v_1,\dots,v_t,u_{t+1},\dots,u_{q'}))
\psi(\e(u_1+\mu_1,\dots,u_t+\mu_t,v_{t+1}+\mu_{t+1},\dots,
v_{q'}+\mu_{q'}))-\\
-\frac{\theta(v_{q'}-u_{q'}+\mu)}{\theta(v_{q'}-u_{q'})}
\psi(\e(v_1,\dots,v_{q'}))\psi(\e(u_1+\mu_1,\dots,u_{q'}+\mu_{q'}))=
\\
=\sum_{\alpha,\beta\in\Z/n\Z}x_{1-\alpha}x_{1-\beta}\times\\
\times\biggl(
\frac{\theta(v_1-u_1+\mu)}{\theta(v_1-u_1)}
w_\alpha(u_1,\dots,u_{q'})
w_\beta(v_1+\mu_1,\dots,
v_{q'}+\mu_{q'})-\\
-\sum_{1\le t<q'}\frac{\theta(\mu)\theta(v_t-u_t+u_{t+1}-v_{t+1})}
{\theta(v_t-u_t)\theta(u_{t+1}-v_{t+1})}\times\\
\times w_\alpha(v_1,\dots,v_t,u_{t+1},\dots,u_{q'})
w_\beta(u_1+\mu_1,\dots,u_t+\mu_t,v_{t+1}+\mu_{t+1},\dots,
v_{q'}+\mu_{q'})-\\
-\frac{\theta(v_{q'}-u_{q'}+\mu)}{\theta(v_{q'}-u_{q'})}
w_\alpha(v_1,\dots,v_{q'})w_\beta(u_1+\mu_1,\dots,u_{q'}+\mu_{q'})
\biggr).
\end{multline*}
By using the identity (35) in Appendix B together with the relations (18) in the 
algebra $Q_{n,k}(\E,\eta)$, one can readily see that this expression is equal to $0$.
\end{proof}

\newpage

\appendix{Theta functions of one variable}

Let $\Gamma\subset\C$ be an integral lattice generated by $1$ and
$\tau\in\C$, where $\Im\tau>0$. Let $n\in\N$ and $c\in\C$. We denote
by  $\Theta_{n,c}(\Gamma)$ the space of the entire functions
of one variable satisfying the following relations:
$$
f(z+1)=f(z),\quad f(z+\tau)=(-1)^ne^{-2\pi i(nz-c)}f(z).
$$
As is known [30], $\dim\Theta_{n,c}(\Gamma)=n$, every
function $f\in\Theta_{n,c}(\Gamma)$ has exactly $n$ zeros modulo $\Gamma$   
(counted according to their multiplicities), and the sum of these zeros 
 modulo $\Gamma$ is equal
to $c$. Let  $\theta(z)=\sum_{\alpha\in\Z}(-1)^\alpha
e^{2\pi i\left(\alpha
z+\frac{\alpha(\alpha-1)}2\tau\right)}$. It is clear that
$\theta(z)\in\Theta_{1,0}(\Gamma)$. It follows from what was
 said above that  $\theta(0)=0$, and this is the only zero 
modulo~$\Gamma$. One can readily see that
$\theta(-z)=-e^{-2\pi iz}\theta(z)$. Moreover, as is known, the function
$\theta(z)$ can be expanded as the infinite product as follows:
$$
\theta(z)=\prod_{\alpha\ge1}(1-e^{2\pi
i\alpha\tau})\cdot(1-e^{2\pi
iz})\cdot\prod_{\alpha\ge1}(1-e^{2\pi
i(z+\alpha\tau)})(1-e^{2\pi i(\alpha\tau-z)}).
$$
 
Let us introduce the following linear operators $T_{\frac1n}$ and $T_{\frac1n\tau}$ 
acting on the space of functions of one variable:
$$
T_{\frac1n}f(z)=f\left(z+\frac1n\right),\quad
T_{\frac1n\tau}f(z)=e^{2\pi
i\left(z+\frac1{2n}-\frac{n-1}{2n}\tau\right)}
f\left(z+\frac1n\tau\right).
$$
One can readily see that the space $\Theta_{n,\frac{n-1}2}(\Gamma)$ is invariant 
with respect to the operators $T_{\frac1n}$ and $T_{\frac1n\tau}$. Moreover, 
$T_{\frac1n}T_{\frac1n\tau}=e^{\frac{2\pi i}n}T_{\frac1n\tau}T_{\frac1n}$. The 
restriction of these operators to the space $\Theta_{n,\frac{n-1}2}(\Gamma)$ satisfy 
the relations $T_{\frac1n}^n=T_{\frac1n\tau}^n=1$. Let $\wt{\Gamma_n}$ be the 
group with the generators $a,b,\epsilon$ and the defining relations 
$ab=\epsilon  ba$, $a\epsilon=\epsilon  a$,  $b\epsilon=\epsilon  b$, and
$a^n=b^n=\epsilon^n=e$. The group $\wt{\Gamma_n}$ is a central extension of the 
group $\Gamma_n=\Gamma/n\Gamma\simeq(\Z/n\Z)^2$, namely, the element $\epsilon$ 
generates a normal subgroup $C_n=\Z/n\Z$, and $\wt{\Gamma_n}/C_n=\Gamma_n$. The 
formulas $a\mapsto   T_{\frac1n}$,
$b\mapsto   T_{\frac1n\tau}$, and  $\epsilon\mapsto{}$(multiplication by
$e^{\frac{2\pi   i}n}$) define an irreducible representation of the group 
$\wt{\Gamma_n}$ in the space $\Theta_{n,\frac{n-1}2}(\Gamma)$. Let us choose a 
basis $\{\theta_\alpha;\alpha\in\Z/n\Z\}$ in the space $\Theta_{n,\frac{n-1}2}(\Gamma)$ 
in which our operators act as follows: $T_{\frac1n}\theta_\alpha=e^{2\pi i\frac\alpha
n}\theta_\alpha$, and $T_{\frac1n\tau}\theta_\alpha=\theta_{\alpha+1}$. It is 
clear that this choice can be carried out uniquely up to multiplication by a 
common constant. The functions $\theta_\alpha(z)$ are of the form 
$$
\theta_\alpha(z)=\theta\left(z+\frac\alpha n\tau\right)
\theta\left(z+\frac1n+\frac\alpha n\tau\right)\dots
\theta\left(z+\frac{n-1}n+\frac\alpha n\tau\right)
e^{2\pi                                                 i\left(\alpha
z+\frac{\alpha(\alpha-n)}{2n}\tau+\frac\alpha{2n}\right)}.
$$
One can readily see that $\theta_\alpha(z)\in\Theta_{n,\frac{n-1}2}(\Gamma)$,
$\theta_{\alpha+n}(z)=\theta_\alpha(z)$, and
\begin{equation}
\begin{aligned}
\theta_\alpha\left(z+\frac1n\right)&=e^{2\pi            i\frac\alpha
n}\theta_\alpha(z),\\
\theta_\alpha\left(z+\frac1n\tau\right)&=e^{-2\pi
i\left(z+\frac1{2n}-\frac{n-1}{2n}\tau\right)}\theta_{\alpha+1}(z)
\end{aligned}
\end{equation}
It is clear that the functions $\left\{\theta_\alpha\left(z-\frac1nc-\frac{n-1}{2n}\right);
\alpha\in\Z/n\Z\right\}$ form a basis in the space $\Theta_{n,c}(\Gamma)$.

We need some identities:
\begin{equation}
\theta(nz)=\frac{n\theta_0(z)\dots\theta_{n-1}(z)e^{-2\pi
i\frac{n(n-1)}2z}}{\theta_1(0)\dots\theta_{n-1}(0)
\theta\left(\frac1n\right)\dots\theta\left(\frac{n-1}n\right)}.
\end{equation}

\begin{proof}
One can readily see by using relations (26) that the functions on both sides of 
the equation belong to the space $\Theta_{n^2,\frac{n(n-1)}2\tau}(\Gamma)$. 
Moreover, it is clear that the zeros of both functions coincide, namely, these 
are $n^2$ points $\left\{\frac\alpha
n+\frac\beta  n\tau;\alpha,\beta\in\Z\right\}$ modulo $\Gamma$. Hence, the 
functions on the left- and right-hand sides of the equation differ by a constant 
multiple, which can be evaluated by dividing (27) by $\theta(z)$ and passing to the 
limit as $z\to0$.
\end{proof}

Let $\theta_0,\theta_1,\theta_2\in\Theta_{3,0}(\Gamma)$. For $z,\eta\in\C$ and 
$\alpha\in\Z/3\Z$ we have
\begin{equation}
\theta_0(\eta)\theta_\alpha(z+\eta)\theta_\alpha(z)+
\theta_1(\eta)\theta_{\alpha+2}(z+\eta)
\theta_{\alpha+1}(z)+
\theta_2(\eta)\theta_{\alpha+1}(z+\eta)
\theta_{\alpha+2}(z)=0.
\end{equation}

\begin{proof}
It is clear that $\theta_\alpha(z+\eta)\theta_\beta(z)\in\Theta_{6,-3\eta}(\Gamma)$ 
as a function of the variable $z$. There must be three linear relations among 
these nine functions in a six-dimensional space. With regard to the action of the 
group $\wt{\Gamma_3}$, we see that the relations must be of the form
$a(\eta)\theta_\alpha(z+\eta)\theta_\alpha(z)+
b(\eta)\theta_{\alpha+1}(z+\eta)\theta_{\alpha+2}(z)+
c(\eta)\theta_{\alpha+2}(z+\eta)\theta_{\alpha+1}(z)=0$. Really, every  
three-dimensional space of relations invariant with respect to the translations 
$z\to  z+\frac13$ and $z\to z+\frac13\tau$ (see (26)) is of this form, where 
$a,b,c$ do not depend on $\alpha$. By setting  $\alpha=1$ and $z=0$, we obtain 
$\frac{c(\eta)}{a(\eta)}=\frac{\theta_1(\eta)}{\theta_0(\eta)}$. By setting 
$\alpha=2$ and $z=0$, we obtain $\frac{b(\eta)}{a(\eta)}=\frac{\theta_2(\eta)}{\theta_0(\eta)}$.
\end{proof}  

Let $\theta_\alpha\in\Theta_{n,c}(\Gamma)$. Then
\begin{multline}
\frac{\theta(y-z+nv-nu)}{\theta(y-z)\theta(nv-nu)}
\theta_\alpha(y+u)\theta_\beta(z+v+\eta)+
\frac{\theta(z-y+n\eta)}{\theta(z-y)\theta(n\eta)}
\theta_\alpha(z+u)\theta_\beta(y+v+\eta)=\\
=\frac1n\theta\left(\frac1n\right)\dots\theta\left(\frac{n-1}n\right)
\sum_{r\in\Z/n\Z}\frac{\theta_{\beta-\alpha}(v-u+\eta)}
{\theta_r(\eta)\theta_{\beta-\alpha-r}(v-u)}\theta_{\beta-r}(y+v)
\theta_{\alpha+r}(z+u+\eta).
\end{multline}

\begin{proof}
This is a special case of the relation (31) (for $p=1$) proved in Appendix B.
\end{proof}

By setting $u=v+\eta$ in the relation (29) and making the change of variables 
$y+v\to  y$, $z+v\to z$, we obtain
\begin{multline}
\frac{\theta(z-y+n\eta)}{\theta(z-y)\theta(n\eta)}
(\theta_\alpha(z+\eta)\theta_\beta(y+\eta)-\theta_\alpha(y+\eta)
\theta_\beta(z+\eta))=\\
=\frac1n\theta\left(\frac1n\right)\dots\theta\left(\frac{n-1}n\right)
\theta_{\beta-\alpha}(0)\sum_{r\in\Z/n\Z}
\frac1{\theta_r(\eta)\theta_{\beta-\alpha-r}(-\eta)}
\theta_{\beta-r}(y)\theta_{\alpha+r}(z+2\eta).
\end{multline}

\appendix{Some theta functions of several variables associated with a power of an 
elliptic curve}

Let $n$ and $k$ be coprime positive integers such
that $1\le
k<n$. We expand the ratio $\frac nk$ in a continued fraction of the
form: $\frac
nk=n_1-\frac1{n_2-\frac1{n_3-\ldots-\frac1{n_p}}}$,
where $n_\alpha\ge2$ for any $\alpha$. It is clear that such an 
expansion exists and is unique. We denote by $d(m_1, \dots , m_q)$ the determinant of the ($q
\times q$) matrix~$(m_{\alpha\beta})$, where
$m_{\alpha\alpha}=m_\alpha$,
$m_{\alpha,\alpha+1}=m_{\alpha+1,\alpha}=-1$, and
$m_{\alpha,\beta}=0$ for $|\alpha-\beta|>1$.  For 
$q = 0$ we set $d(\varnothing)=1$. It follows from
the elementary theory of continued fractions that $n = d(n_1, \dots ,
n_p)$ and $k = d(n_2, \dots , n_p)$.

       Let $\Gamma\subset\C$ be an integral lattice   
generated by 1 and~$\tau$ again, where $\Im\tau>0$.

      We denote by $\Theta_{n/k}(\Gamma)$ the space of entire functions (of $p$ 
variables) satisfying the following relations:
\begin{align*}
f(z_1,\dots,z_\alpha+1,\dots,z_p)&=f(z_1,\dots,z_p),\\
f(z_1,\dots,z_\alpha+\tau,\dots,z_p)&=(-1)^{n_\alpha}e^{-2\pi
i(n_\alpha
z_\alpha-z_{\alpha-1}-z_{\alpha+1}-(\delta_{1,\alpha}-1)\tau)}
f(z_1,\dots,z_p).
\end{align*}
Here $1\le\alpha\le p$ and $z_0=z_{p+1}=0$, and $\delta_{1,\alpha}$
stands for the Kronecker delta. Thus, the functions $f\in\Theta_{n/k}(\Gamma)$
are periodic with respect to each of the variables with period 1 and quasiperiodic with period~$\tau$. 
By the periodicity, each  function in the space $\Theta_{n/k}(\Gamma)$ can be expanded 
in a Fourier series of the form
$f(z_1,\dots,z_p)=\sum_{\alpha_1,\dots,\alpha_p\in\Z}
a_{\alpha_1\dots\alpha_p}e^{2\pi i(\alpha_1z_1+\ldots+\alpha_pz_p)}$. 
By the quasiperiodicity, the coefficients satisfy the system of linear 
equations
$$
a_{\alpha_1,\dots,\alpha_{\nu-1}-1,\alpha_\nu+n_\nu,
\alpha_{\nu+1}-1,\dots,\alpha_p}=
(-1)^{n_\alpha}e^{2\pi i(\alpha_\nu+\delta_{1,\alpha}-1)\tau}
a_{\alpha_1\dots\alpha_p}.
$$

       This system clearly has $n = d(n_1, \dots ,
n_p)$ linearly independent solutions
each defining (for $\Im\tau>0$;
$n_1,\dots,n_p\ge2$) a function in the space $\Theta_{n/k}(\Gamma)$ .

      For $k = 1$ we have the space
of functions of one variable $\Theta_n(\Gamma)=\Theta_{n,0}(\Gamma)$ (see Appendix A) with
a basis
$\left\{w_\alpha(z)=\theta_\alpha\left(z+\frac{n-1}2\right),\
\alpha\in\Z/n\Z\right\}$. A similar basis can be constructed in the space
$\Theta_{n/k}(\Gamma)$ for an arbitrary~$k$. Let us define the operators
$T_{\frac1n}$ and~$T_{\frac1n\tau}$ in the space of functions of $p$
variables as follows:
\begin{align*}
T_{\frac1n}f(z_1,\dots,z_p)&=f(z_1+r_1,\dots,z_p+r_p),\\
T_{\frac1n\tau}f(z_1,\dots,z_p)&=e^{2\pi i(z_1+\phi)}
f(z_1+r_1\tau,\dots,z_p+r_p\tau).
\end{align*}
Here $r_\alpha=\frac{d(n_{\alpha+1},\dots,n_p)}{d(n_1,\dots,n_p)}$ and   
$\phi\in\C$ is a constant.

    It is clear that $T_{\frac1n}T_{\frac1n\tau}=e^{2\pi i\frac
kn}T_{\frac1n\tau}T_{\frac1n}$. As in the case of theta functions of one  
variable, the space
 $\Theta_{n/k}(\Gamma)$ is invariant with respect to the operators $T_{\frac1n}$
and $T_{\frac1n\tau}$, and the restriction of these operators
to~$\Theta_{n/k}(\Gamma)$ satisfies the relations $T_{\frac1n}^n=1$ and 
$T_{\frac1n\tau}^n=\mu$, where $\mu\in\C$. Let us choose a $\phi$  in such a way that
$\mu=1$; clearly, this can be done uniquely up to multiplication of
$T_{\frac1n\tau}$ by a root of unity of degree $n$.

\begin{prop} 
There is a basis $\bigl\{w_\alpha(z_1,\dots,z_p);\,\alpha\in\Z/n\Z\bigr\}$ in  
$\Theta_{n/k}(\Gamma)$ such that
$$
T_{\frac1n}w_\alpha=e^{2\pi i\frac kn\alpha}w_\alpha,\quad
T_{\frac1n\tau}w_\alpha=w_{\alpha+1}.
$$
This basis is defined uniquely up to multiplication by a common constant.
\end{prop}

\begin{proof}
 Let $f\in\Theta_{n/k}(\Gamma)$ be an eigenvector of the operator $T_{\frac1n}$ with 
 an eigenvalue~$\lambda$. Since $T_{\frac1n}^n=1$ on the space $\Theta_{n/k}(\Gamma)$,
we have  $\lambda^n=1$. Moreover, $T_{\frac1n}T_{\frac1n\tau}f=e^{2\pi
i\frac
kn}T_{\frac1n\tau}T_{\frac1n}f=e^{2\pi i\frac kn}\lambda
T_{\frac1n\tau}f$, and hence $T_{\frac1n\tau}f$ is also
an eigenvector with the eigenvalue $e^{2\pi i\frac
kn}\lambda$. Since $n$ and $k$ are coprime, $e^{2\pi i\frac
kn}$ is a
primitive root of unity of degree $n$. Thus, the vectors $\bigl\{T_{\frac1n\tau}^\alpha
f;\,\alpha=0,1,\dots,n-1\bigr\}$
 are eigenvectors for the operator $T_{\frac1n}$ with
different eigenvalues, and every of root of unity of degree $n$ is an
eigenvalue for some $T_{\frac1n\tau}^\alpha f$. Let $w_0$
be such that $T_{\frac1n}w_0=w_0$. We set
$w_\alpha=T_{\frac1n\tau}^\alpha w_0$. It is clear that
$T_{\frac1n}w_\alpha=e^{2\pi i\frac
kn\alpha}w_\alpha$ and
$T_{\frac1n\tau}w_\alpha=w_{\alpha+1}$. Moreover,
$w_{\alpha+n}=w_\alpha$ because $T_{\frac1n\tau}^n=1$
on the space $\Theta_{n/k}(\Gamma)$.
\end{proof}

    We note that, as in the case of theta functions of one variable, the group 
    $\wt{\Gamma_n}$ irreducibly acts on
the space $\Theta_{n/k}(\Gamma)$ by the rule $a\mapsto T_{\frac1n}$, $b\mapsto
T_{\frac1n\tau}$, and $\epsilon\mapsto\text{(multiplication by $e^{2\pi
i\frac kn}$})$.

\begin{remark}
 Let $L$ be the group of linear automorphisms on the space of functions 
 of $p$
variables of the form
$$
gf(z_1,\dots,z_p)=e^{2\pi
i(\phi_1z_1+\ldots+\phi_pz_p+\lambda)}f(z_1+\psi_1,\dots,z_p+\psi_p)
$$
for $g\in L$. It is clear that $L$ is a $(2p + 1)$-dimensional Lie
group. Let $L'\subset L$ be the subgroup of
transformations preserving the space~$\Theta_{n/k}(\Gamma)$, that is,
$L'=\bigl\{g\in
L;\,g(\Theta_{n/k}(\Gamma))=\Theta_{n/k}(\Gamma)\bigr\}$.
Let $L''\subset L'$  consist of the elements preserving each point of
$\Theta_{n/k}(\Gamma)$, that is, $L''=\bigl\{g\in
L';\,gf=f\text{ for any }f\in\Theta_{n/k}(\Gamma)\bigr\}$. One can see that 
the quotient group $L'/L''=\wt G_n$  
 is generated by the elements
$T_{\frac1n}$ and $T_{\frac1n\tau}$ and by the multiplications by
constants.
\end{remark}
We shall use the notation $w_\alpha^{n/k}(z_1,\dots,z_p)$ if it is not clear 
from the context what are the theta functions in use.

We need the following identity relating theta functions in the spaces 
$\Theta_{1,0}(\Gamma)$, $\Theta_{n,\frac{n-1}2}(\Gamma)$, and $\Theta_{n/k}(\Gamma)$: 
\begin{align}
&\frac{\theta(y_1-z_1+nv-nu)}{\theta(nv-nu)\theta(y_1-z_1)}\nonumber\\
&{}\times w_\alpha(y_1+m_1u,\dots,y_p+m_pu)
w_\beta(z_1+m_1v+l_1,\dots,z_p+m_pv+l_p)\nonumber\\
&+\sum_{1\le t\le p-1}
\frac{\theta(z_t-y_t+y_{t+1}-z_{t+1})}
{\theta(z_t-y_t)\theta(y_{t+1}-z_{t+1})}\nonumber\\
&{}\times w_\alpha(z_1+m_1u,\dots,z_t+m_tu,y_{t+1}+m_{t+1}u,\dots,y_p+m_pu)
\nonumber\\ &{}\times
w_\beta(y_1+m_1v+l_1,\dots,y_t+m_tv+l_t,z_{t+1}+m_{t+1}v+l_{t+1},
\dots,z_p+m_pv+l_p)\nonumber\\
&{}+\frac{\theta(z_p-y_p+n\eta)}{\theta(z_p-y_p)\theta(n\eta)}
w_\alpha(z_1+m_1u_1,\dots,z_p+m_pu)
\nonumber\\ &{}\times w_\beta(y_1+m_1v+l_1,\dots,y_p+m_pv+l_p)\nonumber\\
&{}=\frac1n\theta\left(\frac1n\right)\dots\theta\left(\frac{n-1}n\right)
\nonumber\\ &{}\times
\sum_{r\in\Z/n\Z}\frac{\theta_{\beta-\alpha+ r(k-1)}(v-u+\eta)}
{\theta_{ r k}(\eta)\theta_{\beta-\alpha- r}(v-u)}
w_{\beta- r}(y_1+m_1v,\dots,y_p+m_pv)
\nonumber\\ &{}\times
w_{\alpha+ r}(z_1+m_1u+l_1,\dots,z_p+m_pu+l_p).
\end{align}
Here $m_\alpha=d(n_{\alpha+1},\dots,n_p)$  and 
$l_\alpha=d(n_1,\dots,n_{\alpha-1})\eta$.

\begin{proof}
 We denote by   
$\phi_{\alpha,\beta}(\eta,u,v,y_1,\dots,y_p,z_1,\dots,z_p)$
the difference between the right- and left-hand 
sides of the formula (31). The calculation shows that this function 
satisfies the following relations:
\begin{equation}
\begin{aligned}
\phi_{\alpha,\beta}(\eta,\dots,y_\alpha+1,\dots,z_p)&=
\phi_{\alpha,\beta}(\eta,\dots,z_p),\\
\phi_{\alpha,\beta}(\eta,\dots,y_\alpha+\tau,\dots,z_p)&=
-e^{-2\pi i(n_\alpha y_\alpha-y_{\alpha-1}-y_{\alpha+1}+\delta_{\alpha,1}v)}
\phi_{\alpha,\beta}(\eta,\dots,z_p),\\
\phi_{\alpha,\beta}(\eta,\dots,z_\alpha+1,\dots,z_p)&=
\phi_{\alpha,\beta}(\eta,\dots,z_p),\\
\phi_{\alpha,\beta}(\eta,\dots,z_\alpha+\tau,\dots,z_p)&=
-e^{-2\pi i(n_\alpha z_\alpha-z_{\alpha-1}-z_{\alpha+1}+
\delta_{\alpha,1}u+\delta_{\alpha,p}\eta)}
\phi_{\alpha,\beta}(\eta,\dots,z_p).
\end{aligned}
\end{equation}
Here $y_0=y_{p+1}=z_0=z_{p+1}=0$, and   $\delta_{\alpha,\beta}$  stands for the Kronecker
 delta. Moreover, an evaluation shows
that there are no poles on the divisors $nv-nu\in\Gamma$,
$n\eta\in\Gamma$,
$y_1-z_1\in\Gamma$, \dots, $y_p-z_p\in\Gamma$,
and hence the function   $\phi_{\alpha,\beta}$   is holomorphic everywhere 
on~$\C^{2p+3}$. However, it is clear that the
functions $\bigl\{w_\lambda(y_1+m_1v,\dots,y_p+m_pv)
w_\nu(z_1+m_1u+l_1,\dots,z_p+m_pu+l_p);\,
\lambda,\nu\in\Z/n\Z\bigr\}$ form
a basis in the space of holomorphic functions (of the variables $y_1,\dots,y_p,z_1,\dots,z_p$) 
 satisfying the conditions (32).
Therefore, the function   $\phi_{\alpha,\beta}$  is of the form
\begin{multline}
\varphi_{\alpha,\beta}(\eta,u,v,y_1,\dots,z_p)=\\
{}=\sum_{\lambda,\nu\in\Z/n\Z}
\psi_{\lambda,\nu}(\eta,u,v)w_\lambda(y_1+m_1v,\dots,y_p+m_pv)
\\ {}\times
w_\nu(z_1+m_1u+l_1,\dots,z_p+m_pu+l_p).
\end{multline}

        Here the functions   $\psi_{\lambda,\nu}(\eta,u,v)$ are holomorphic and 
satisfy the relations
\begin{equation}
\begin{gathered}
\psi_{\lambda,\nu}(\eta+1,u,v)=\psi_{\lambda,\nu}(\eta,u+1,v)=
\psi_{\lambda,\nu}(\eta,u,v+1)=\psi_{\lambda,\nu}(\eta,u,v),\\
\psi_{\lambda,\nu}(\eta+\tau,u,v)=
e^{-2\pi in(v-u)}\psi_{\lambda,\nu}(\eta,u,v),\\
\psi_{\lambda,\nu}(\eta,u+\tau,v)=e^{2\pi in\eta}\psi_{\lambda,\nu}(\eta,u,v),\\
\psi_{\lambda,\nu}(\eta,u,v+\tau)=e^{-2\pi in\eta}\psi_{\lambda,\nu}(\eta,u,v).
\end{gathered}
\end{equation}

        These relations are verified by the immediate calculation, namely, 
one must compare the multipliers at the translations by 1
and~$\tau$   in the formulas (32) and (33).

        However, every holomorphic function of the variables $\eta$, $u$ and $v$
 that satisfies relations (34) is vanishes. 
 Really, since this function is periodic, it admits the expansion in the Fourier 
series
$$
\psi_{\lambda,\nu}(\eta,u,v)=\sum_{\alpha,\beta,\gamma\in\Z}
a_{\lambda,\nu,\alpha,\beta,\gamma}
e^{2\pi i(\alpha\eta+\beta u+\gamma v)}.
$$
Further, it follows from the quasiperiodicity that the
coefficients $a_{\lambda,\nu,\alpha,\beta,\gamma}$  are equal to 0.
\end{proof}

By setting $u=v+\eta$ in the identity (31) and making the change of variables 
$y_1\to
y_1-m_1v$,  $z_1\to  z_1-m_1v$,  \dots,  $y_p\to  y_p-m_pv$,  $z_p\to
z_p-m_pv$, we obtain
\begin{multline}
\frac{\theta(y_1-z_1-n\eta)}{\theta(-n\eta)\theta(y_1-z_1)}
w_\alpha(y_1+m_1\eta,\dots,y_p+m_p\eta)
\theta_\beta(z_1+\l_1,\dots,z_p+\l_p)+\\
+\sum_{1\le t<p}\frac{\theta(z_t-y_t+y_{t+1}-z_{t+1})}
{\theta(z_t-y_t)\theta(y_{t+1}-z_{t+1})}
w_\alpha(z_1+m_1\eta,\dots,z_t+m_t\eta,
y_{t+1}+m_{t+1}\eta,\dots,y_p+m_p\eta)\times\\
\times
w_\beta(y_1+\l_1,\dots,y_t+\l_t,z_{t+1}+\l_{t+1},\dots,z_p+\l_p)+\\
+\frac{\theta(z_p-y_p+n\eta)}{\theta(z_p-y_p)\theta(n\eta)}
w_\alpha(z_1+m_1\eta,\dots,z_p+m_p\eta)
w_\beta(y_1+\l_1,\dots,y_p+\l_p)=\\
=\frac1n\theta\left(\frac1n\right)\dots\theta\left(\frac{n-1}n\right)
\times\\\times\sum_{r\in\Z/n\Z}\frac{\theta_{\beta-\alpha+r(k-1)}(0)}
{\theta_{rk}(\eta)\theta_{\beta-\alpha-r}(-\eta)}
w_{\beta-r}(y_1,\dots,y_p)
w_{\alpha+r}(z_1+m_1\eta+\l_1,\dots,z_p+m_p\eta+\l_p).
\end{multline}

\appendix{Duality between the spaces  $\Theta_{n/k}(\Gamma)$ and
$\Theta_{n/n-k}(\Gamma)$} 

Let us construct a canonical element 
$\Delta_{n,k}\in\Theta_{n/k}(\Gamma)\otimes\Theta_{n/n-k}(\Gamma)$
carring out the duality between these spaces (see (36)).

\begin{prop}
 Let
$$
\frac nk=n_1-\frac1{n_2-\ldots-\frac1{n_p}},\quad \frac
n{n-k}=n_1'-\frac1{n_2'-\ldots-\frac1{n_{p'}'}}
$$ 
be the expansions in continued fractions, where $n_\alpha\ge2$ and
$n_\beta'\ge2$  for $1\le\alpha\le p$ and $1\le\beta\le p'$, respectively.
Here $p$ and $p'$  stand for the lengths of the continued fractions. Then 
$p'=n_1+\ldots+n_p-2p+1$
and $n_1'+\ldots+n_{p'}'=2(n_1+\ldots+n_p)-3p+1$. Moreover,
$n_1'+\ldots+n_\alpha'=2\alpha+\beta$   for
$n_1+\ldots+n_\beta-2\beta+1\le\alpha\le
n_1+\ldots+n_{\beta+1}-2\beta-2$. In other words, the Young diagrams for the 
partitions
$(n_1-1,n_1+n_2-3,\dots,n_1+\ldots+n_\alpha-2\alpha+1,\dots)$ and 
$(n_1'-1,n_1'+n_2'-3,\dots,n_1'+\dots+n_\beta'-2\beta+1,\dots)$
are dual to each other.
\end{prop}

\begin{remark} 
For $k = 1$, $p = 1$, and 
$n_1 = n$ we have $p'  = n - 1$ and $n_1'=\ldots=n_{n-1}'=2$. 
For $p > 1$, if $n_2, \dots , n_{p-1}\ge   3$, then the sequence
$(n_1',\dots,n_p')$ becomes
$(2^{(n_1-2)},3,2^{(n_2-3)},3,\dots,3,2^{(n_{p-1}-3)},3,2^{(n_p-2)})$.
Here $2^{(t)}$, $t\ge0$, stands for a sequence of $t$ twos. This
formula remains valid without the assumption that $n_2,\dots,n_p\ge3$ if we agree that
 the
sequence $(m_1,2^{(-1)},m_2)$ is of length $1$ and is equal to $(m_1+m_2-2)$. 
This rule must be
applied in succession to all occurrences $n_\alpha=2$ for $2\le\alpha\le p-1$.
\end{remark}

\begin{proof}
 The proof can be carried out by induction on $\min(p,p')$. For $p = 1$, one must 
 prove that
$\frac n{n-1}=2-\frac1{2-\ldots-\frac12}$  is of length 
$n - 1$. Let $p,p'>1$ and let, say, $n_1>2$. We
have     $\frac
k{d(n_3,\dots,n_p)}=n_2-\frac1{n_3-\ldots-\frac1{n_p}}$. By assumption,
$$
\frac k{k-d(n_3,\dots,n_p)}=
n_{n_1-1}'-1-\frac1{n_{n_1}'-\frac1{n_{n_1+1}'-\ldots-\frac1{n_{p'}'}}}.
$$
and the sequence $(n_1',\dots,n_{n_1-2}')$ is~$(2^{(n_1-2)})$.
Further, one must show that
$n_1'-\frac1{n_2'-\ldots-\frac1{n_{p'}'}}=\frac n{n-k}$. Here it is 
used that 
$d(n_1,\dots,n_p)=n$,
$d(n_2,\dots,n_p)=k$,
$
\frac{d(n_1,\dots,n_p)}{d(n_2,\dots,n_p)}=
n_1-\frac1{n_2-\ldots-\frac1{n_p}},
$ 
and $d(n_1,\dots,n_p)=n_1d(n_2,\dots,n_p)-d(n_3,\dots,n_p)$.
\end{proof}

\begin{prop} 
Let a function  $\Delta_{n.k}(z_1,\dots,z_p;z_1',\dots,z_{p'}')$
 of $p+p'$  variables $z_1,\dots,z_p,z_1',\dots,z_{p'}'$
  be defined by the formula
\begin{multline*}
\Delta_{n,k}(z_1,\dots,z_p,z_1',\dots,z_{p'}')=\\
{}=e^{2\pi
iz_1'}\theta(z_1-z_1')\theta(z_p+z_{p'}')\cdot\prod_{1\le\alpha\le
p'-1}\theta(z_\alpha'-z_{\alpha+1}'+z_{n_1'+\ldots+n_\alpha'-2\alpha+1})
\\
{}\times\prod_{1\le\beta\le
p-1}\theta(z_\beta-z_{\beta+1}+z_{n_1+\ldots+n_\beta-2\beta+1}').
\end{multline*}
This function satisfies the following relations:
$$
\Delta_{n,k}(z_1,\dots,z_\alpha+1,\dots,z_{p'}')=
\Delta_{n,k}(z_1,\dots,z_\beta'+1,\dots,z_{p'}')=
\Delta_{n,k}(z_1,\dots,z_{p'}'),\\
$$
\begin{multline*}
\Delta_{n,k}(z_1,\dots,z_\alpha+\tau,\dots,z_{p'}')\\
=(-1)^{n_\alpha}e^{-2\pi i
(n_\alpha z_\alpha-z_{\alpha-1}-z_{\alpha+1}-(\delta_{\alpha,1}-1)\tau)}
\Delta_{n,k}(z_1,\dots,z_{p'}'),
\end{multline*}
\begin{multline*}
\Delta_{n,k}(z_1,\dots,z_\beta'+\tau,\dots,z_{p'}')\\
=(-1)^{n_\beta'}e^{-2\pi i
(n_\beta'z_\beta'-z_{\beta-1}'-z_{\beta+1}'-(\delta_{\beta,1}-1)\tau)}
\Delta_{n,k}(z_1,\dots,z_{p'}').
\end{multline*}
Here $z_0=z_{p+1}=z_0'=z_{p'+1}'=0$ and $\delta_{\alpha,1}$  stands for the
Kronecker delta.
\end{prop}

      The proof immediately follows from our description of the duality 
 between the sequences
$(n_1,\dots,n_p)$ and $(n_1',\dots,n_{p'}')$.

\begin{prop}                                                           
\begin{equation}
\Delta_{n,k}(z_1,\dots,z_p;z_1',\dots,z_{p'}')=c_{n,k}
\sum_{\alpha\in\Z/n\Z}w_\alpha^{n/k}(z_1,\dots,z_p)
w_{1-\alpha}^{n/n-k}(z_1',\dots,z_{p'}').
\end{equation}
Here $c_{n,k}\in\C$ is a constant.
\end{prop}

\begin{proof} 
It follows from the previous proposition that the function $\Delta_{n,k}$  
 belongs to the space~$\Theta_{n/k}(\Gamma)$ when regarded as a function of the 
variables $z_1,\dots,z_p$. Similarly, $\Delta_{n,k}$  belongs to~$\Theta_{n/n-k}(\Gamma)$
as a function of $z_1',\dots,z_{p'}'$. Therefore,  
$$
\Delta_{n,k}(z_1,\dots,z_p;z_1',\dots,z_{p'}')=
\sum_{\alpha,\beta\in\Z/n\Z}\lambda_{\alpha,\beta}
w_\alpha^{n/k}(z_1,\dots,z_p)
w_\beta^{n/n-k}(z_1',\dots,z_{p'}').
$$
However, one can readily see that
$$
\Delta_{n,k}(z_1+ r_1,\dots,z_p+ r_p;z_1'+ r_1',\dots,z_{p'}'+ r_{p'}')=
e^{\frac{2\pi i}n}\Delta_{n,k}(z_1,\dots,z_{p'}'),
$$
where $ r_\alpha=\frac{d(n_1,\dots,n_{\alpha-1})}n$ and $
r_\beta'=\frac{d(n_1',\dots,n_{\beta-1}')}n$. Hence,
$\lambda_{\alpha,\beta}=0$ for  $\alpha+\beta\not\equiv1\mmod n$
(because $w_\alpha(z_1+ r_1,\dots,z_p+ r_p)=e^{2\pi i\frac\alpha
n}w_\alpha(z_1,\dots,z_p)$ and $w_\beta(z_1'+ r_1',\dots,z_{p'}'+
r_{p'}')=e^{2\pi i\frac\beta
n}w_\beta(z_1',\dots,z_{p'}')$). Thus,   
$\lambda_{\alpha,\beta}=\lambda_\alpha\delta_{\alpha+\beta,1}$. Similarly,
\begin{multline*}
\Delta_{n,k}(z_1+ r_1\tau,\dots,z_p+ r_p\tau;z_1'+ r_1'\tau,
\dots,z_{p'}'+ r_{p'}'\tau)\\
{}=e^{2\pi i\left(\frac1n\tau-z_p-z_{p'}'\right)}
\Delta_{n,k}(z_1,\dots,z_{p'}').
\end{multline*}
Hence, $\lambda_\alpha=\lambda_{\alpha+1}$, that is,
$\lambda_\alpha$    does not depend on~$\alpha$.
\end{proof}

\appendix{}

\subsection{Integrable system, quantum groups, and $R$-matrices}

One of the main methods in the investigation of exactly solvable models [6] in 
quantum and statistical physics is the inverse problem method (see [45]). This 
method leads to the study of representations of algebras of monodromy matrices, 
that is, to the study of meromorphic matrix functions $L(u)$ satisfying the relations
\begin{equation}
R(u-v)L^1(u)L^2(v)=L^2(v)L^1(u)R(u-v).
\end{equation}
Here $R(u)$ is a chosen solution of the Yang-Baxter equation in the class of 
meromorphic matrix-valued functions,
\begin{equation}
R^{12}(u-v)R^{13}(u)R^{23}(v)=R^{23}(v)R^{13}(u)R^{12}(u-v).
\end{equation}
We note that $R(u)$ takes the values in ($n^2\times n^2$) matrices with a fixed 
decomposition $\Mat_{n^2}=\Mat_n\otimes\Mat_n$. We use the standard notation, 
namely, $L^1=L\otimes1$, $L^2=1\otimes L$, $R^{12}=R\otimes1$, etc. (see [45]). 

In [45] Sklyanin studies the solutions of the equation (37) for the simplest 
elliptic solution of the equation (38), that is, for the so-called Baxter $R$-matrix, 
which is of the form
 $R(u)=
1+\sum_{\alpha=1}^3W_\alpha(u)\sigma_\alpha\otimes\sigma_\alpha$, where
$\sigma_1=\begin{pmatrix}0&1\\1&0\end{pmatrix}$,
$\sigma_2=\begin{pmatrix}0&-i\\i&0\end{pmatrix}$, and
$\sigma_3=\begin{pmatrix}1&0\\0&-1\end{pmatrix}$ are the Pauli matrices, and the 
coefficients $W_\alpha(u)$ can be expressed in terms of the Jacobi elliptic 
functions as follows:
\begin{align*}
W_1(u)&=\frac{\sn(i\eta,k)}{\sn(u+i\eta,k)},\quad
W_2(u)=\frac{\dn}{\sn}(u+i\eta,k)\frac{\sn}{\dn}(i\eta,k),\\
W_3(u)&=\frac{\cn}{\sn}(u+i\eta,k)\frac{\sn}{\cn}(i\eta,k).
\end{align*}
The functions $W_\alpha(u)$ uniformize the elliptic curve 
$\frac{W_\alpha^2-W_\beta^2}{W_\gamma^2-1}=J_{\alpha,\beta}$, where the 
$J_{\alpha,\beta}$s do not depend on $u$ and satisfy the relation 
$J_{12}+J_{23}+J_{31}+J_{12}J_{23}J_{31}=0$. Here $\alpha,\beta$, and $\gamma$ 
are pairwise distinct and $J_{\beta,\alpha}=-J_{\alpha,\beta}$. 

We note that this elliptic curve is the complete intersection of two quadrics, 
for instance, $w_1^2-w_2^2=J_{12}(w_3^2-1)$ and $w_2^2-w_3^2=J_{23}(w_1^2-1)$. 

Sklyanin discovered that the equation (37) for the Baxter $R$-matrix has a solution 
of the form $L(u)=S_0+\sum_{\alpha=1}^3W_\alpha(u)S_\alpha$, where $S_0$ 
and $S_\alpha$ are matrices that do not depend on $u$ and satisfy the following 
relations:
\begin{equation}
\begin{aligned}\relax
[S_\alpha,S_0]_-&=-iJ_{\beta\gamma}[S_\beta,S_\gamma]_+,\\
[S_\alpha,S_\beta]_-&=i[S_0,S_\gamma]_+,\\
\end{aligned}
\end{equation}
where $[a,b]_\pm=ab\pm ba$. 

Sklyanin further studies the algebra with the generators $S_0,S_\alpha$ and the 
relations (39); he denotes this algebra by $\F_{\eta,k}$. The main assumption 
concerning this algebra is that it satisfies the PBW condition. Moreover, Sklyanin 
finds the quadratic central elements of the algebra $\F_{\eta,k}$ and the 
finite-dimensional representations of the algebra $\F_{\eta,k}$ by difference 
operators in some function space (see [46]).

In our notation, the Sklyanin algebra $\F_{\eta,k}$ is the algebra $Q_4(\E,\eta)$, 
where $\E$ is an elliptic curve given by the functions $W_\alpha(u)$, that is, a 
complete intersection of two quadrics in $\C^3$.

The Yang-Baxter equation has other elliptic solutions generalizing the Baxter 
solution (see [7]). The result of [7] can be described as follows: for any pair 
of positive integers $n$ and $k$ such that $1\le  k<n$ and $n$ and $k$ are coprime 
there is a family of solutions $R_{n,k}(\E,\eta)(u)$ of the equation (38). Here 
$\E$ is an elliptic curve and $\eta\in\E$, as above. The Baxter solution is obtained 
for $n=2$ and $k=1$.

According to [10], the Sklyanin result can be generalized to an arbitrary solution 
$R_{n,k}(\E,\eta)(u)$. In our notation, the results of [10] look as follows: there 
is a homomorphism of the algebra of monodromy matrices for the $R$-matrix 
$R_{n,k}(\E,\eta)$ into the algebra $Q_{n^2,nk-1}(\E,\eta)$. Correspondingly, the 
algebra $Q_{n^2,nk-1}(\E,\eta)$ is a deformation of the projectivization of the 
Lie algebra $\slg_n$. Moreover, there is a homomorphism of the algebra of monodromy 
matrices into the algebra $Q_{dn^2,dnk-1}(\E,\eta)$ for any $d\in\N$. It can be 
conjectured that every finite-dimensional representation of the algebra of monodromy 
matrices can be obtained from a representation of the algebra $Q_{dn^2,dnk-1}(\E,\eta)$. 

Another relationship between the elliptic solutions of the Yang-Baxter equation 
and the elliptic algebras follows from the results of [10]. The multiplication in 
the algebra $Q_{n,k}(\E,\eta)$ is defined by the so-called Young projections 
$S^\alpha V\otimes S^\beta V\to S^{\alpha+\beta}V$ corresponding to 
$R_{n,k}(\E,\eta)(u)$ (see [10]). Moreover, $Q_{n,k}(\E,\eta)=\sum_\alpha S^\alpha V$.

We also note that the study of algebras $Q_{n,k}(\E,\eta)$ and their representations 
led to deeper understanding of the structure of $R$-matrices $R_{n,k}(\E,\eta)(u)$ 
and of the corresponding algebraic objects (the Zamolodchikov algebra and the 
algebra of monodromy matrices). For this topic, see [33].

\subsection{Deformation quantization}

Let $M$ be a manifold ($C^\infty$, analytic, algebraic, etc.) and let $\F(M)$ 
be a function algebra on $M$. In the "physical" language, $M$ is the state space 
of the system and $\F(M)$ is the algebra of observables. In [5] the following 
approach to the quantization was suggested: the underlying vector space of the 
quantum algebra of observables coincide with that of $\F(M)$, but the multiplication 
is deformed and is no longer commutative (though still associative). Moreover, 
the multiplication depends on the deformation parameter (Planck constant). For 
$\hbar=0$ we have the ordinary commutative multiplication. Since the Planck 
constant is small, we do not notice that the observables in classical mechanics 
are non-commutative. Expanding the multiplication in the series in powers of $\hbar$ 
we obtain $f*g=fg+\{f,g\}\hbar+o(\hbar)$. The operation 
$\{{\cdot},{\cdot}\}\colon\F(M)\otimes\F(M)\to\F(M)$ is bilinear, and, applying 
the gauge transformations, one can make it anticommutative, $\{f,g\}=-\{g,f\}$. 
Moreover, since the multiplication $*$ is associative, we see that 
$\{f,gh\}=\{fg\}h+\{f,h\}g$ (the Leibniz rule) and 
$\{f,\{g,h\}\}+\{h,\{f,g\}\}+\{g,\{h,f\}\}=0$ (the Jacobi identity). The Leibniz 
rule means that $\{f,g\}=\langle w,df\wedge dg\rangle$, where $w$ is a bivector 
field on $M$, and the Jacobi identity means that $[w,w]=0$. Thus, every quantization 
defines a Poisson structure $w$ (or $\{{\cdot},{\cdot}\}$) on the manifold $M$. 
The inverse problem arises: construct a quantization $*$ from a manifold $M$ 
with Poisson structure (a Poisson manifold). This problem was solved in [26] at 
the formal level. Namely, bidifferential operators $B_n$  ($n\ge2$) were constructed 
on a Poisson manifold $M$ in such a way that the formal series 
\begin{equation}
f*g=fg+\{f,g\}\hbar+\sum_{n\ge2}B_n(f,g)\hbar^n
\end{equation}
satisfies the condition $(f*g)*h=f*(g*h)$. Moreover, the operators $B_n$ are 
constructed from $\{{\cdot},{\cdot}\}$ by an explicit formula. Thus, the problem 
of formal quantization was solved; however, the problem on the convergence of the 
series (40) and on its identification remains open. This problem seems to be very 
complicated. For instance, let $M=\C^n$,  $\F(M)=S^*(\C^n)=\C[x_1,\dots,x_n]$, 
and let the Poisson bracket be quadratic ($\{x_i,x_j\}=\sum_{\alpha,\beta}c_{ij}^{\alpha\beta}x_\alpha
x_\beta$, where $c_{ij}^{\alpha\beta}\in\C$ are symmetric with respect to $\alpha,\beta$ 
and antisymmetric with respect to $i,j$). In this case, the multiplication $*$ 
must be homogeneous, that is, $S^\alpha(\C^n)*S^\beta(\C^n)\subset
S^{\alpha+\beta}(\C^n)$. Therefore, according to (40), the structure constants 
of the multiplication $*$ in the basis $\{x_1^{\alpha_1}\dots
x_n^{\alpha_n};\alpha_1,\dots,\alpha_n\in\Z_{\ge0}\}$ turn out to be formal series 
in $\hbar$ which baffle the explicit evaluation even in the simplest case $n=2$,
$\{x_1,x_2\}=\alpha x_1x_2$. In this case it is natural to assume that the quantum 
algebra must be defined by the relation $x_1x_2=e^{-\alpha\hbar}x_2x_1$. On the 
other hand, the algebras $Q_{n,k}(\E,\eta)$ introduced above are examples of the 
quantization of $M=\C^n$, where $\eta$ plays the role of Planck constant because 
$Q_{n,k}(\E,0)=\C[x_1,\dots,x_n]$. Moreover, the structure constants of the 
algebra $Q_{n,k}(\E,\eta)$ turn out to be elliptic functions of $\eta$. There are 
also rational and trigonometric limits of the algebras $Q_{n,k}(\E,\eta)$ in which 
the structure constants are rational (trigonometric) functions of $\eta$ (see [32]).

\subsection{Moduli spaces}

Let $G$ be a semisimple Lie group and let $P\subset G$ be a parabolic subgroup. 
Let $\M(\E,P)$ be the moduli space of the holomorphic $P$-bundles over an elliptic 
curve $\E$ [4]. According to [20], every connected component of the space $\M(\E,P)$ 
admits a natural Poisson structure. The main property of this structure is as 
follows: the preimages of the natural map $\M(\E,P)\to\M(\E,G)$ corresponding to 
forgetting the $P$-structure are symplectic leaves of the structure. The quantization 
problem for the Poisson manifold $\M(\E,P)$ arises. The solution of this problem 
could establish a relationship between the natural algebro-geometric problem of 
studying $P$-bundles (and the corresponding $G$-bundles) and the problem to study 
representations of the quantum function algebra on $\M(\E,P)$ because the representations 
correspond to symplectic leaves.

In [38], [22] these quantum algebras were constructed provided that $P=B$ is a 
Borel subgroup of an arbitrary group $G$. In [21] the quantum algebras were 
constructed in the case of $G=GL_m$ and an arbitrary parabolic subgroup $P$. Here 
the algebra $Q_n(\E,\eta)$ corresponds to the case $G=GL_2$, and the algebra 
$Q_{n,k}(\E,\eta)$ to the case $G=GL_{k+1}$, where $P$ consists of upper block 
triangular matrices of the form
$\vcenter{\raisebox{5pt}{\bm{\mathstrut\cr
k&&*&&*\cr\mathstrut\cr1&&0&&*\cr&&k&&1}}
\vspace{20pt}.}$

\subsection{Non-commutative algebraic geometry}

One of the main ideas of algebraic geometry is to study the geometry of a manifold 
by using the algebraic properties of a ring of functions on this manifold. The 
non-commutative algebraic geometry extends these methods and the geometric intuition 
to an appropriate class of non-commutative rings. In [51], [48] the non-commutative 
algebraic geometry is developed in small dimensions. From this point of view, 
the algebras $Q_{n,k}(\E,\eta)$ give examples of non-commutative vector spaces. 
Similar examples of non-commutative Grassmannians ans other varieties are also 
known [11], [19], [22], [49].

\subsection{Cohomology of algebras}

Cohomology properties of quadratic algebras are studied in [29], [41]-[44], [50]. 
For generic $\eta$, the algebras $Q_{n,k}(\E,\eta)$ are examples of Koszul algebras. 
One can readily prove this fact for $k=1$ by using the construction of a free 
module in \S2.6. The constructions of dual algebras $Q_{n,k}^!(\E,\eta)$ are 
given in [36].

\bigskip
\begin{flushleft}
L. D. Landau Institute of Theoretical Physics,\\
Russian Academy of Sciences\\
\medskip
e-mail: odesskii@mccme.ru

\end{flushleft}

\typeout{LaTeX Warning: Label(s) may have changed. Rerun}
\end{document}